\input amstex
\input amsppt.sty
\magnification1200

\def\LL{\leavevmode\setbox0=\hbox{L}\hbox to\wd0{\hss\char'40L}}
\def\al{\alpha}

\def\ze{\zeta}

\def\la{\lambda}
\def\rh{\rho}

\def\om{\omega}


\def\today{\ifcase\month\or
 January\or February\or March\or April\or May\or June\or
 July\or August\or September\or October\or November\or December\fi
 \space\number\day, \number\year}

\def\({\left(}
\def\){\right)}
\def\[{\left[}
\def\]{\right]}

\def\3{\ss}
\catcode`\@=11
\def\dddot#1{\vbox{\ialign{##\crcr
      .\hskip-.5pt.\hskip-.5pt.\crcr\noalign{\kern1.5\p@\nointerlineskip}
      $\hfil\displaystyle{#1}\hfil$\crcr}}}

\newif\iftab@\tab@false
\newif\ifvtab@\vtab@false
\def\tab{\bgroup\tab@true\vtab@false\vst@bfalse\Strich@false%
   \def\\{\global\hline@@false%
     \ifhline@\global\hline@false\global\hline@@true\fi\cr}
   \edef\l@{\the\leftskip}\ialign\bgroup\hskip\l@##\hfil&&##\hfil\cr}
\def\endtab{\cr\egroup\egroup}
\def\vtab{\vtop\bgroup\vst@bfalse\vtab@true\tab@true\Strich@false%
   \bgroup\def\\{\cr}\ialign\bgroup&##\hfil\cr}
\def\endvtab{\cr\egroup\egroup\egroup}
\def\stab{\D@cke0.5pt\null 
 \bgroup\tab@true\vtab@false\vst@bfalse\Strich@true\Let@@\vspace@
 \normalbaselines\offinterlineskip
  \openup\spreadmlines@
 \edef\l@{\the\leftskip}\ialign
 \bgroup\hskip\l@##\hfil&&##\hfil\crcr}
\def\endstab{\crcr\egroup
 \egroup}
\newif\ifvst@b\vst@bfalse
\def\vstab{\D@cke0.5pt\null
 \vtop\bgroup\tab@true\vtab@false\vst@btrue\Strich@true\bgroup\Let@@\vspace@
 \normalbaselines\offinterlineskip
  \openup\spreadmlines@\bgroup}
\def\endvstab{\crcr\egroup\egroup
 \egroup\tab@false\Strich@false}

\newdimen\htstrut@
\htstrut@8.5\p@
\newdimen\htStrut@
\htStrut@12\p@
\newdimen\dpstrut@
\dpstrut@3.5\p@
\newdimen\dpStrut@
\dpStrut@3.5\p@
\def\openup{\afterassignment\@penup\dimen@=}
\def\@penup{\advance\lineskip\dimen@
  \advance\baselineskip\dimen@
  \advance\lineskiplimit\dimen@
  \divide\dimen@ by2
  \advance\htstrut@\dimen@
  \advance\htStrut@\dimen@
  \advance\dpstrut@\dimen@
  \advance\dpStrut@\dimen@}
\def\Let@@{\relax%
    \def\\{\global\hline@@false%
     \ifhline@\global\hline@false\global\hline@@true\fi\cr}%
    \iffalse}\fi}
\def\matrix{\null\,\vcenter\bgroup
 \tab@false\vtab@false\vst@bfalse\Strich@false\Let@@\vspace@
 \normalbaselines\openup\spreadmlines@\ialign
 \bgroup\hfil$\m@th##$\hfil&&\quad\hfil$\m@th##$\hfil\crcr
 \Mathstrut@\crcr\noalign{\kern-\baselineskip}}
\def\endmatrix{\crcr\Mathstrut@\crcr\noalign{\kern-\baselineskip}\egroup
 \egroup\,}
\def\smatrix{\D@cke0.5pt\null\,
 \vcenter\bgroup\tab@false\vtab@false\vst@bfalse\Strich@true\Let@@\vspace@
 \normalbaselines\offinterlineskip
  \openup\spreadmlines@\ialign
 \bgroup\hfil$\m@th##$\hfil&&\quad\hfil$\m@th##$\hfil\crcr}
\def\endsmatrix{\crcr\egroup
 \egroup\,\Strich@false}
\newdimen\D@cke
\def\Dicke#1{\global\D@cke#1}
\newtoks\tabs@\tabs@{&}
\newif\ifStrich@\Strich@false
\newif\iff@rst

\def\Stricherr@{\iftab@\ifvtab@\errmessage{\noexpand\s not allowed
     here. Use \noexpand\vstab!}%
  \else\errmessage{\noexpand\s not allowed here. Use \noexpand\stab!}%
  \fi\else\errmessage{\noexpand\s not allowed
     here. Use \noexpand\smatrix!}\fi}
\def\format{\ifvst@b\else\crcr\fi\egroup\iffalse{\fi\ifnum`}=0 \fi\format@}
\def\format@#1\\{\def\preamble@{#1}%
 \def\Str@chfehlt##1{\ifx##1\s\Stricherr@\fi\ifx##1\\\let\Next\relax%
   \else\let\Next\Str@chfehlt\fi\Next}%
 \def\c{\hfil\noexpand\ifhline@@\hbox{\vrule height\htStrut@%
   depth\dpstrut@ width\z@}\noexpand\fi%
   \ifStrich@\hbox{\vrule height\htstrut@ depth\dpstrut@ width\z@}%
   \fi\iftab@\else$\m@th\fi\the\hashtoks@\iftab@\else$\fi\hfil}%
 \def\r{\hfil\noexpand\ifhline@@\hbox{\vrule height\htStrut@%
   depth\dpstrut@ width\z@}\noexpand\fi%
   \ifStrich@\hbox{\vrule height\htstrut@ depth\dpstrut@ width\z@}%
   \fi\iftab@\else$\m@th\fi\the\hashtoks@\iftab@\else$\fi}%
 \def\l{\noexpand\ifhline@@\hbox{\vrule height\htStrut@%
   depth\dpstrut@ width\z@}\noexpand\fi%
   \ifStrich@\hbox{\vrule height\htstrut@ depth\dpstrut@ width\z@}%
   \fi\iftab@\else$\m@th\fi\the\hashtoks@\iftab@\else$\fi\hfil}%
 \def\s{\ifStrich@\ \the\tabs@\vrule width\D@cke\the\hashtoks@%
          \fi\the\tabs@\ }%
 \def\sa{\ifStrich@\vrule width\D@cke\the\hashtoks@%
            \the\tabs@\ %
            \fi}%
 \def\se{\ifStrich@\ \the\tabs@\vrule width\D@cke\the\hashtoks@\fi}%
 \def\cd{\hfil\noexpand\ifhline@@\hbox{\vrule height\htStrut@%
   depth\dpstrut@ width\z@}\noexpand\fi%
   \ifStrich@\hbox{\vrule height\htstrut@ depth\dpstrut@ width\z@}%
   \fi$\dsize\m@th\the\hashtoks@$\hfil}%
 \def\rd{\hfil\noexpand\ifhline@@\hbox{\vrule height\htStrut@%
   depth\dpstrut@ width\z@}\noexpand\fi%
   \ifStrich@\hbox{\vrule height\htstrut@ depth\dpstrut@ width\z@}%
   \fi$\dsize\m@th\the\hashtoks@$}%
 \def\ld{\noexpand\ifhline@@\hbox{\vrule height\htStrut@%
   depth\dpstrut@ width\z@}\noexpand\fi%
   \ifStrich@\hbox{\vrule height\htstrut@ depth\dpstrut@ width\z@}%
   \fi$\dsize\m@th\the\hashtoks@$\hfil}%
 \ifStrich@\else\Str@chfehlt#1\\\fi%
 \setbox\z@\hbox{\xdef\Preamble@{\preamble@}}\ifnum`{=0 \fi\iffalse}\fi
 \ialign\bgroup\span\Preamble@\crcr}
\newif\ifhline@\hline@false
\newif\ifhline@@\hline@@false
\def\hlinefor#1{\multispan@{\strip@#1 }\leaders\hrule height\D@cke\hfill%
    \global\hline@true\ignorespaces}
\def\Item "#1"{\par\noindent\hangindent2\parindent%
  \hangafter1\setbox0\hbox{\rm#1\enspace}\ifdim\wd0>2\parindent%
  \box0\else\hbox to 2\parindent{\rm#1\hfil}\fi\ignorespaces}
\def\ITEM #1"#2"{\par\noindent\hangafter1\hangindent#1%
  \setbox0\hbox{\rm#2\enspace}\ifdim\wd0>#1%
  \box0\else\hbox to 0pt{\rm#2\hss}\hskip#1\fi\ignorespaces}
\def\item"#1"{\par\noindent\hang%
  \setbox0=\hbox{\rm#1\enspace}\ifdim\wd0>\the\parindent%
  \box0\else\hbox to \parindent{\rm#1\hfil}\enspace\fi\ignorespaces}
\let\plainitem@\item
\catcode`\@=13

\catcode`\@=11
\font@\twelverm=cmr10 scaled\magstep1
\font@\twelveit=cmti10 scaled\magstep1
\font@\twelvebf=cmbx10 scaled\magstep1
\font@\twelvei=cmmi10 scaled\magstep1
\font@\twelvesy=cmsy10 scaled\magstep1
\font@\twelveex=cmex10 scaled\magstep1

\newtoks\twelvepoint@
\def\twelvepoint{\normalbaselineskip15\p@
 \abovedisplayskip15\p@ plus3.6\p@ minus10.8\p@
 \belowdisplayskip\abovedisplayskip
 \abovedisplayshortskip\z@ plus3.6\p@
 \belowdisplayshortskip8.4\p@ plus3.6\p@ minus4.8\p@
 \textonlyfont@\rm\twelverm \textonlyfont@\it\twelveit
 \textonlyfont@\sl\twelvesl \textonlyfont@\bf\twelvebf
 \textonlyfont@\smc\twelvesmc \textonlyfont@\tt\twelvett
%
 \ifsyntax@ \def\big##1{{\hbox{$\left##1\right.$}}}%
  \let\Big\big \let\bigg\big \let\Bigg\big
 \else
  \textfont\z@=\twelverm  \scriptfont\z@=\tenrm  \scriptscriptfont\z@=\sevenrm
  \textfont\@ne=\twelvei  \scriptfont\@ne=\teni  \scriptscriptfont\@ne=\seveni
  \textfont\tw@=\twelvesy \scriptfont\tw@=\tensy \scriptscriptfont\tw@=\sevensy
  \textfont\thr@@=\twelveex \scriptfont\thr@@=\tenex
        \scriptscriptfont\thr@@=\tenex
  \textfont\itfam=\twelveit \scriptfont\itfam=\tenit
        \scriptscriptfont\itfam=\tenit
  \textfont\bffam=\twelvebf \scriptfont\bffam=\tenbf
        \scriptscriptfont\bffam=\sevenbf
  \setbox\strutbox\hbox{\vrule height10.2\p@ depth4.2\p@ width\z@}%
  \setbox\strutbox@\hbox{\lower.6\normallineskiplimit\vbox{%
        \kern-\normallineskiplimit\copy\strutbox}}%
 \setbox\z@\vbox{\hbox{$($}\kern\z@}\bigsize@=1.4\ht\z@
 \fi
 \normalbaselines\rm\ex@.2326ex\jot3.6\ex@\the\twelvepoint@}

\font@\fourteenrm=cmr10 scaled\magstep2
\font@\fourteenit=cmti10 scaled\magstep2
\font@\fourteensl=cmsl10 scaled\magstep2
\font@\fourteensmc=cmcsc10 scaled\magstep2
\font@\fourteentt=cmtt10 scaled\magstep2
\font@\fourteenbf=cmbx10 scaled\magstep2
\font@\fourteeni=cmmi10 scaled\magstep2
\font@\fourteensy=cmsy10 scaled\magstep2
\font@\fourteenex=cmex10 scaled\magstep2
\font@\fourteenmsa=msam10 scaled\magstep2
\font@\fourteeneufm=eufm10 scaled\magstep2
\font@\fourteenmsb=msbm10 scaled\magstep2
\newtoks\fourteenpoint@
\def\fourteenpoint{\normalbaselineskip15\p@
 \abovedisplayskip18\p@ plus4.3\p@ minus12.9\p@
 \belowdisplayskip\abovedisplayskip
 \abovedisplayshortskip\z@ plus4.3\p@
 \belowdisplayshortskip10.1\p@ plus4.3\p@ minus5.8\p@
 \textonlyfont@\rm\fourteenrm \textonlyfont@\it\fourteenit
 \textonlyfont@\sl\fourteensl \textonlyfont@\bf\fourteenbf
 \textonlyfont@\smc\fourteensmc \textonlyfont@\tt\fourteentt
%
 \ifsyntax@ \def\big##1{{\hbox{$\left##1\right.$}}}%
  \let\Big\big \let\bigg\big \let\Bigg\big
 \else
  \textfont\z@=\fourteenrm  \scriptfont\z@=\twelverm  \scriptscriptfont\z@=\tenrm
  \textfont\@ne=\fourteeni  \scriptfont\@ne=\twelvei  \scriptscriptfont\@ne=\teni
  \textfont\tw@=\fourteensy \scriptfont\tw@=\twelvesy \scriptscriptfont\tw@=\tensy
  \textfont\thr@@=\fourteenex \scriptfont\thr@@=\twelveex
        \scriptscriptfont\thr@@=\twelveex
  \textfont\itfam=\fourteenit \scriptfont\itfam=\twelveit
        \scriptscriptfont\itfam=\twelveit
  \textfont\bffam=\fourteenbf \scriptfont\bffam=\twelvebf
        \scriptscriptfont\bffam=\tenbf
  \setbox\strutbox\hbox{\vrule height12.2\p@ depth5\p@ width\z@}%
  \setbox\strutbox@\hbox{\lower.72\normallineskiplimit\vbox{%
        \kern-\normallineskiplimit\copy\strutbox}}%
 \setbox\z@\vbox{\hbox{$($}\kern\z@}\bigsize@=1.7\ht\z@
 \fi
 \normalbaselines\rm\ex@.2326ex\jot4.3\ex@\the\fourteenpoint@}

\font@\seventeenrm=cmr10 scaled\magstep3
\font@\seventeenit=cmti10 scaled\magstep3
\font@\seventeensl=cmsl10 scaled\magstep3
\font@\seventeensmc=cmcsc10 scaled\magstep3
\font@\seventeentt=cmtt10 scaled\magstep3
\font@\seventeenbf=cmbx10 scaled\magstep3
\font@\seventeeni=cmmi10 scaled\magstep3
\font@\seventeensy=cmsy10 scaled\magstep3
\font@\seventeenex=cmex10 scaled\magstep3
\font@\seventeenmsa=msam10 scaled\magstep3
\font@\seventeeneufm=eufm10 scaled\magstep3
\font@\seventeenmsb=msbm10 scaled\magstep3
\newtoks\seventeenpoint@
\def\seventeenpoint{\normalbaselineskip18\p@
 \abovedisplayskip21.6\p@ plus5.2\p@ minus15.4\p@
 \belowdisplayskip\abovedisplayskip
 \abovedisplayshortskip\z@ plus5.2\p@
 \belowdisplayshortskip12.1\p@ plus5.2\p@ minus7\p@
 \textonlyfont@\rm\seventeenrm \textonlyfont@\it\seventeenit
 \textonlyfont@\sl\seventeensl \textonlyfont@\bf\seventeenbf
 \textonlyfont@\smc\seventeensmc \textonlyfont@\tt\seventeentt
%
 \ifsyntax@ \def\big##1{{\hbox{$\left##1\right.$}}}%
  \let\Big\big \let\bigg\big \let\Bigg\big
 \else
  \textfont\z@=\seventeenrm  \scriptfont\z@=\fourteenrm  \scriptscriptfont\z@=\twelverm
  \textfont\@ne=\seventeeni  \scriptfont\@ne=\fourteeni  \scriptscriptfont\@ne=\twelvei
  \textfont\tw@=\seventeensy \scriptfont\tw@=\fourteensy \scriptscriptfont\tw@=\twelvesy
  \textfont\thr@@=\seventeenex \scriptfont\thr@@=\fourteenex
        \scriptscriptfont\thr@@=\fourteenex
  \textfont\itfam=\seventeenit \scriptfont\itfam=\fourteenit
        \scriptscriptfont\itfam=\fourteenit
  \textfont\bffam=\seventeenbf \scriptfont\bffam=\fourteenbf
        \scriptscriptfont\bffam=\twelvebf
  \setbox\strutbox\hbox{\vrule height14.6\p@ depth6\p@ width\z@}%
  \setbox\strutbox@\hbox{\lower.86\normallineskiplimit\vbox{%
        \kern-\normallineskiplimit\copy\strutbox}}%
 \setbox\z@\vbox{\hbox{$($}\kern\z@}\bigsize@=2\ht\z@
 \fi
 \normalbaselines\rm\ex@.2326ex\jot5.2\ex@\the\seventeenpoint@}

\catcode`\@=13
\catcode`\@=11
\font\tenln    = line10
\font\tenlnw   = linew10

\newskip\Einheit \Einheit=0.5cm
\newcount\xcoord \newcount\ycoord
\newdimen\xdim \newdimen\ydim \newdimen\PfadD@cke \newdimen\Pfadd@cke

\newcount\@tempcnta
\newcount\@tempcntb

\newdimen\@tempdima
\newdimen\@tempdimb

\newdimen\@wholewidth
\newdimen\@halfwidth

\newcount\@xarg
\newcount\@yarg
\newcount\@yyarg
\newbox\@linechar
\newbox\@tempboxa
\newdimen\@linelen
\newdimen\@clnwd
\newdimen\@clnht

\newif\if@negarg

\def\@whilenoop#1{}
\def\@whiledim#1\do #2{\ifdim #1\relax#2\@iwhiledim{#1\relax#2}\fi}
\def\@iwhiledim#1{\ifdim #1\let\@nextwhile=\@iwhiledim
        \else\let\@nextwhile=\@whilenoop\fi\@nextwhile{#1}}

\def\@whileswnoop#1\fi{}
\def\@whilesw#1\fi#2{#1#2\@iwhilesw{#1#2}\fi\fi}
\def\@iwhilesw#1\fi{#1\let\@nextwhile=\@iwhilesw
         \else\let\@nextwhile=\@whileswnoop\fi\@nextwhile{#1}\fi}

\def\thinlines{\let\@linefnt\tenln \let\@circlefnt\tencirc
  \@wholewidth\fontdimen8\tenln \@halfwidth .5\@wholewidth}
\def\thicklines{\let\@linefnt\tenlnw \let\@circlefnt\tencircw
  \@wholewidth\fontdimen8\tenlnw \@halfwidth .5\@wholewidth}
\thinlines

\PfadD@cke1pt \Pfadd@cke0.5pt
\def\PfadDicke#1{\PfadD@cke#1 \divide\PfadD@cke by2 \Pfadd@cke\PfadD@cke \multiply\PfadD@cke by2}
\long\def\LOOP#1\REPEAT{\def\BODY{#1}\ITERATE}
\def\ITERATE{\BODY \let\next\ITERATE \else\let\next\relax\fi \next}
\let\REPEAT=\fi
\def\Punkt{\hbox{\raise-2pt\hbox to0pt{\hss$\ssize\bullet$\hss}}}
\def\DuennPunkt(#1,#2){\unskip
  \raise#2 \Einheit\hbox to0pt{\hskip#1 \Einheit
          \raise-2.5pt\hbox to0pt{\hss$\bullet$\hss}\hss}}
\def\NormalPunkt(#1,#2){\unskip
  \raise#2 \Einheit\hbox to0pt{\hskip#1 \Einheit
          \raise-3pt\hbox to0pt{\hss\twelvepoint$\bullet$\hss}\hss}}
\def\DickPunkt(#1,#2){\unskip
  \raise#2 \Einheit\hbox to0pt{\hskip#1 \Einheit
          \raise-4pt\hbox to0pt{\hss\fourteenpoint$\bullet$\hss}\hss}}
\def\Kreis(#1,#2){\unskip
  \raise#2 \Einheit\hbox to0pt{\hskip#1 \Einheit
          \raise-4pt\hbox to0pt{\hss\fourteenpoint$\circ$\hss}\hss}}

\def\Line@(#1,#2)#3{\@xarg #1\relax \@yarg #2\relax
\@linelen=#3\Einheit
\ifnum\@xarg =0 \@vline
  \else \ifnum\@yarg =0 \@hline \else \@sline\fi
\fi}

\def\@sline{\ifnum\@xarg< 0 \@negargtrue \@xarg -\@xarg \@yyarg -\@yarg
  \else \@negargfalse \@yyarg \@yarg \fi
\ifnum \@yyarg >0 \@tempcnta\@yyarg \else \@tempcnta -\@yyarg \fi
\ifnum\@tempcnta>6 \@badlinearg\@tempcnta0 \fi
\ifnum\@xarg>6 \@badlinearg\@xarg 1 \fi
\setbox\@linechar\hbox{\@linefnt\@getlinechar(\@xarg,\@yyarg)}%
\ifnum \@yarg >0 \let\@upordown\raise \@clnht\z@
   \else\let\@upordown\lower \@clnht \ht\@linechar\fi
\@clnwd=\wd\@linechar
\if@negarg \hskip -\wd\@linechar \def\@tempa{\hskip -2\wd\@linechar}\else
     \let\@tempa\relax \fi
\@whiledim \@clnwd <\@linelen \do
  {\@upordown\@clnht\copy\@linechar
   \@tempa
   \advance\@clnht \ht\@linechar
   \advance\@clnwd \wd\@linechar}%
\advance\@clnht -\ht\@linechar
\advance\@clnwd -\wd\@linechar
\@tempdima\@linelen\advance\@tempdima -\@clnwd
\@tempdimb\@tempdima\advance\@tempdimb -\wd\@linechar
\if@negarg \hskip -\@tempdimb \else \hskip \@tempdimb \fi
\multiply\@tempdima \@m
\@tempcnta \@tempdima \@tempdima \wd\@linechar \divide\@tempcnta \@tempdima
\@tempdima \ht\@linechar \multiply\@tempdima \@tempcnta
\divide\@tempdima \@m
\advance\@clnht \@tempdima
\ifdim \@linelen <\wd\@linechar
   \hskip \wd\@linechar
  \else\@upordown\@clnht\copy\@linechar\fi}

\def\@hline{\ifnum \@xarg <0 \hskip -\@linelen \fi
\vrule height\Pfadd@cke width \@linelen depth\Pfadd@cke
\ifnum \@xarg <0 \hskip -\@linelen \fi}

\def\@getlinechar(#1,#2){\@tempcnta#1\relax\multiply\@tempcnta 8
\advance\@tempcnta -9 \ifnum #2>0 \advance\@tempcnta #2\relax\else
\advance\@tempcnta -#2\relax\advance\@tempcnta 64 \fi
\char\@tempcnta}

\def\Vektor(#1,#2)#3(#4,#5){\unskip\leavevmode
  \xcoord#4\relax \ycoord#5\relax
      \raise\ycoord \Einheit\hbox to0pt{\hskip\xcoord \Einheit
         \Vector@(#1,#2){#3}\hss}}

\def\Vector@(#1,#2)#3{\@xarg #1\relax \@yarg #2\relax
\@tempcnta \ifnum\@xarg<0 -\@xarg\else\@xarg\fi
\ifnum\@tempcnta<5\relax
\@linelen=#3\Einheit
\ifnum\@xarg =0 \@vvector
  \else \ifnum\@yarg =0 \@hvector \else \@svector\fi
\fi
\else\@badlinearg\fi}

\def\@hvector{\@hline\hbox to 0pt{\@linefnt
\ifnum \@xarg <0 \@getlarrow(1,0)\hss\else
    \hss\@getrarrow(1,0)\fi}}

\def\@vvector{\ifnum \@yarg <0 \@downvector \else \@upvector \fi}

\def\@svector{\@sline
\@tempcnta\@yarg \ifnum\@tempcnta <0 \@tempcnta=-\@tempcnta\fi
\ifnum\@tempcnta <5
  \hskip -\wd\@linechar
  \@upordown\@clnht \hbox{\@linefnt  \if@negarg
  \@getlarrow(\@xarg,\@yyarg) \else \@getrarrow(\@xarg,\@yyarg) \fi}%
\else\@badlinearg\fi}

\def\@upline{\hbox to \z@{\hskip -.5\Pfadd@cke \vrule width \Pfadd@cke
   height \@linelen depth \z@\hss}}

\def\@downline{\hbox to \z@{\hskip -.5\Pfadd@cke \vrule width \Pfadd@cke
   height \z@ depth \@linelen \hss}}

\def\@upvector{\@upline\setbox\@tempboxa\hbox{\@linefnt\char'66}\raise
     \@linelen \hbox to\z@{\lower \ht\@tempboxa\box\@tempboxa\hss}}

\def\@downvector{\@downline\lower \@linelen
      \hbox to \z@{\@linefnt\char'77\hss}}

\def\@getlarrow(#1,#2){\ifnum #2 =\z@ \@tempcnta='33\else
\@tempcnta=#1\relax\multiply\@tempcnta \sixt@@n \advance\@tempcnta
-9 \@tempcntb=#2\relax\multiply\@tempcntb \tw@
\ifnum \@tempcntb >0 \advance\@tempcnta \@tempcntb\relax
\else\advance\@tempcnta -\@tempcntb\advance\@tempcnta 64
\fi\fi\char\@tempcnta}

\def\@getrarrow(#1,#2){\@tempcntb=#2\relax
\ifnum\@tempcntb < 0 \@tempcntb=-\@tempcntb\relax\fi
\ifcase \@tempcntb\relax \@tempcnta='55 \or
\ifnum #1<3 \@tempcnta=#1\relax\multiply\@tempcnta
24 \advance\@tempcnta -6 \else \ifnum #1=3 \@tempcnta=49
\else\@tempcnta=58 \fi\fi\or
\ifnum #1<3 \@tempcnta=#1\relax\multiply\@tempcnta
24 \advance\@tempcnta -3 \else \@tempcnta=51\fi\or
\@tempcnta=#1\relax\multiply\@tempcnta
\sixt@@n \advance\@tempcnta -\tw@ \else
\@tempcnta=#1\relax\multiply\@tempcnta
\sixt@@n \advance\@tempcnta 7 \fi\ifnum #2<0 \advance\@tempcnta 64 \fi
\char\@tempcnta}

\def\Diagonale(#1,#2)#3{\unskip\leavevmode
  \xcoord#1\relax \ycoord#2\relax
      \raise\ycoord \Einheit\hbox to0pt{\hskip\xcoord \Einheit
         \Line@(1,1){#3}\hss}}
\def\AntiDiagonale(#1,#2)#3{\unskip\leavevmode
  \xcoord#1\relax \ycoord#2\relax 
      \raise\ycoord \Einheit\hbox to0pt{\hskip\xcoord \Einheit
         \Line@(1,-1){#3}\hss}}
\def\Pfad(#1,#2),#3\endPfad{\unskip\leavevmode
  \xcoord#1 \ycoord#2 \thicklines\ZeichnePfad#3\endPfad\thinlines}
\def\ZeichnePfad#1{\ifx#1\endPfad\let\next\relax
  \else\let\next\ZeichnePfad
    \ifnum#1=1
      \raise\ycoord \Einheit\hbox to0pt{\hskip\xcoord \Einheit
         \vrule height\Pfadd@cke width1 \Einheit depth\Pfadd@cke\hss}%
      \advance\xcoord by 1
    \else\ifnum#1=2
      \raise\ycoord \Einheit\hbox to0pt{\hskip\xcoord \Einheit
        \hbox{\hskip-\PfadD@cke\vrule height1 \Einheit width\PfadD@cke depth0pt}\hss}%
      \advance\ycoord by 1
    \else\ifnum#1=3
      \raise\ycoord \Einheit\hbox to0pt{\hskip\xcoord \Einheit
         \Line@(1,1){1}\hss}
      \advance\xcoord by 1
      \advance\ycoord by 1
    \else\ifnum#1=4
      \raise\ycoord \Einheit\hbox to0pt{\hskip\xcoord \Einheit
         \Line@(1,-1){1}\hss}
      \advance\xcoord by 1
      \advance\ycoord by -1
    \fi\fi\fi\fi
  \fi\next}
\def\hSSchritt{\leavevmode\raise-.4pt\hbox to0pt{\hss.\hss}\hskip.2\Einheit
  \raise-.4pt\hbox to0pt{\hss.\hss}\hskip.2\Einheit
  \raise-.4pt\hbox to0pt{\hss.\hss}\hskip.2\Einheit
  \raise-.4pt\hbox to0pt{\hss.\hss}\hskip.2\Einheit
  \raise-.4pt\hbox to0pt{\hss.\hss}\hskip.2\Einheit}
\def\vSSchritt{\vbox{\baselineskip.2\Einheit\lineskiplimit0pt
\hbox{.}\hbox{.}\hbox{.}\hbox{.}\hbox{.}}}
\def\DSSchritt{\leavevmode\raise-.4pt\hbox to0pt{%
  \hbox to0pt{\hss.\hss}\hskip.2\Einheit
  \raise.2\Einheit\hbox to0pt{\hss.\hss}\hskip.2\Einheit
  \raise.4\Einheit\hbox to0pt{\hss.\hss}\hskip.2\Einheit
  \raise.6\Einheit\hbox to0pt{\hss.\hss}\hskip.2\Einheit
  \raise.8\Einheit\hbox to0pt{\hss.\hss}\hss}}
\def\dSSchritt{\leavevmode\raise-.4pt\hbox to0pt{%
  \hbox to0pt{\hss.\hss}\hskip.2\Einheit
  \raise-.2\Einheit\hbox to0pt{\hss.\hss}\hskip.2\Einheit
  \raise-.4\Einheit\hbox to0pt{\hss.\hss}\hskip.2\Einheit
  \raise-.6\Einheit\hbox to0pt{\hss.\hss}\hskip.2\Einheit
  \raise-.8\Einheit\hbox to0pt{\hss.\hss}\hss}}
\def\SPfad(#1,#2),#3\endSPfad{\unskip\leavevmode
  \xcoord#1 \ycoord#2 \ZeichneSPfad#3\endSPfad}
\def\ZeichneSPfad#1{\ifx#1\endSPfad\let\next\relax
  \else\let\next\ZeichneSPfad
    \ifnum#1=1
      \raise\ycoord \Einheit\hbox to0pt{\hskip\xcoord \Einheit
         \hSSchritt\hss}%
      \advance\xcoord by 1
    \else\ifnum#1=2
      \raise\ycoord \Einheit\hbox to0pt{\hskip\xcoord \Einheit
        \hbox{\hskip-2pt \vSSchritt}\hss}%
      \advance\ycoord by 1
    \else\ifnum#1=3
      \raise\ycoord \Einheit\hbox to0pt{\hskip\xcoord \Einheit
         \DSSchritt\hss}
      \advance\xcoord by 1
      \advance\ycoord by 1
    \else\ifnum#1=4
      \raise\ycoord \Einheit\hbox to0pt{\hskip\xcoord \Einheit
         \dSSchritt\hss}
      \advance\xcoord by 1
      \advance\ycoord by -1
    \fi\fi\fi\fi
  \fi\next}
\def\Koordinatenachsen(#1,#2){\unskip
 \hbox to0pt{\hskip-.5pt\vrule height#2 \Einheit width.5pt depth1 \Einheit}%
 \hbox to0pt{\hskip-1 \Einheit \xcoord#1 \advance\xcoord by1
    \vrule height0.25pt width\xcoord \Einheit depth0.25pt\hss}}
\def\Koordinatenachsen(#1,#2)(#3,#4){\unskip
 \hbox to0pt{\hskip-.5pt \ycoord-#4 \advance\ycoord by1
    \vrule height#2 \Einheit width.5pt depth\ycoord \Einheit}%
 \hbox to0pt{\hskip-1 \Einheit \hskip#3\Einheit 
    \xcoord#1 \advance\xcoord by1 \advance\xcoord by-#3 
    \vrule height0.25pt width\xcoord \Einheit depth0.25pt\hss}}
\def\Gitter(#1,#2){\unskip \xcoord0 \ycoord0 \leavevmode
  \LOOP\ifnum\ycoord<#2
    \loop\ifnum\xcoord<#1
      \raise\ycoord \Einheit\hbox to0pt{\hskip\xcoord \Einheit\Punkt\hss}%
      \advance\xcoord by1
    \repeat
    \xcoord0
    \advance\ycoord by1
  \REPEAT}
\def\Gitter(#1,#2)(#3,#4){\unskip \xcoord#3 \ycoord#4 \leavevmode
  \LOOP\ifnum\ycoord<#2
    \loop\ifnum\xcoord<#1
      \raise\ycoord \Einheit\hbox to0pt{\hskip\xcoord \Einheit\Punkt\hss}%
      \advance\xcoord by1
    \repeat
    \xcoord#3
    \advance\ycoord by1
  \REPEAT}
\def\Label#1#2(#3,#4){\unskip \xdim#3 \Einheit \ydim#4 \Einheit
  \def\lo{\advance\xdim by-.5 \Einheit \advance\ydim by.5 \Einheit}%
  \def\llo{\advance\xdim by-.25cm \advance\ydim by.5 \Einheit}%
  \def\loo{\advance\xdim by-.5 \Einheit \advance\ydim by.25cm}%
  \def\o{\advance\ydim by.25cm}%
  \def\ro{\advance\xdim by.5 \Einheit \advance\ydim by.5 \Einheit}%
  \def\rro{\advance\xdim by.25cm \advance\ydim by.5 \Einheit}%
  \def\roo{\advance\xdim by.5 \Einheit \advance\ydim by.25cm}%
  \def\l{\advance\xdim by-.30cm}%
  \def\r{\advance\xdim by.30cm}%
  \def\lu{\advance\xdim by-.5 \Einheit \advance\ydim by-.6 \Einheit}%
  \def\llu{\advance\xdim by-.25cm \advance\ydim by-.6 \Einheit}%
  \def\luu{\advance\xdim by-.5 \Einheit \advance\ydim by-.30cm}%
  \def\u{\advance\ydim by-.30cm}%
  \def\ru{\advance\xdim by.5 \Einheit \advance\ydim by-.6 \Einheit}%
  \def\rru{\advance\xdim by.25cm \advance\ydim by-.6 \Einheit}%
  \def\ruu{\advance\xdim by.5 \Einheit \advance\ydim by-.30cm}%
  #1\raise\ydim\hbox to0pt{\hskip\xdim
     \vbox to0pt{\vss\hbox to0pt{\hss$#2$\hss}\vss}\hss}%
}
\catcode`\@=13

\NoBlackBoxes
\NoRunningHeads
\TagsOnRight
\loadbold

\def\CoMNAA{1}
\def\FrRTAA{2}
\def\GaMiAB{3}
\def\HiGrAA{4}
\def\KratAY{5}
\def\LuRSAA{6}
\def\MacMAA{7}
\def\NoPSAA{8}
\def\PaStAA{9}
\def\PropAD{10}
\def\PrWiAA{11}
\def\PrWiAB{12}
\def\ReWhAA{13}
\def\RobbAA{14}
\def\SagaAL{15}
\def\SchuAA{16}
\def\StanAA{17}
\def\WilDAA{18}
\def\WilDAC{19}

\def\bP{{\bar P}}
\def\bT{{\bar T}}
\def\bH{{\bar H}}

\def\bom{{\bar\om}}
\def\hC{{C_R}}
\def\hT{{T_L}}
\def\hH{{H_L}}

\catcode`\@=11
\newif\if@ovt
\newif\if@ovb
\newif\if@ovl
\newif\if@ovr
\newdimen\@ovxx
\newdimen\@ovyy
\newdimen\@ovdx
\newdimen\@ovdy
\newdimen\@ovro
\newdimen\@ovri
\newcount\@tempcnta
\newcount\@tempcntb
\newdimen\@tempdima
\newdimen\@tempdimb
\newdimen\@wholewidth
\newdimen\@halfwidth
\newdimen\unitlength \unitlength =1pt

\def\thinlines{\let\@linefnt\tenln \let\@circlefnt\tencirc
  \@wholewidth\fontdimen8\tenln \@halfwidth .5\@wholewidth}
\def\thicklines{\let\@linefnt\tenlnw \let\@circlefnt\tencircw
  \@wholewidth\fontdimen8\tenlnw \@halfwidth .5\@wholewidth}

\font\tenln    = line10
\font\tenlnw   = linew10
\font\tencirc  = lcircle10    
\font\tencircw = lcirclew10   
\thinlines

\def\@ifstar#1#2{\@ifnextchar *{\def\@tempa*{#1}\@tempa}{#2}}
\def\@ifnextchar#1#2#3{\let\@tempe #1\def\@tempa{#2}\def\@tempb{#3}\futurelet
    \@tempc\@ifnch}
\def\@ifnch{\ifx \@tempc \@sptoken \let\@tempd\@xifnch
      \else \ifx \@tempc \@tempe\let\@tempd\@tempa\else\let\@tempd\@tempb\fi
      \fi \@tempd}

\def\@getcirc#1{\@tempdima #1\relax \advance\@tempdima 2pt\relax
  \@tempcnta\@tempdima
  \@tempdima 4pt\relax \divide\@tempcnta\@tempdima
  \ifnum \@tempcnta > 10\relax \@tempcnta 10\relax\fi
  \ifnum \@tempcnta >\z@ \advance\@tempcnta\m@ne
    \else \@warning{Oval too small}\fi
  \multiply\@tempcnta 4\relax
  \setbox \@tempboxa \hbox{\@circlefnt
  \char \@tempcnta}\@tempdima \wd \@tempboxa}

\def\@put#1#2#3{\raise #2\hbox to \z@{\hskip #1#3\hss}}

\def\circle{\@ifstar{\@dot}{\@circle}}
\def\@circle#1{\begingroup \boxmaxdepth \maxdimen \@tempdimb #1\unitlength
   \ifdim \@tempdimb >15.5pt\relax \@getcirc\@tempdimb
      \@ovro\ht\@tempboxa
     \setbox\@tempboxa\hbox{\@circlefnt
      \advance\@tempcnta\tw@ \char \@tempcnta
      \advance\@tempcnta\m@ne \char \@tempcnta \kern -2\@tempdima
      \advance\@tempcnta\tw@
      \raise \@tempdima \hbox{\char\@tempcnta}\raise \@tempdima
        \box\@tempboxa}\ht\@tempboxa\z@ \dp\@tempboxa\z@
      \@put{-\@ovro}{-\@ovro}{\box\@tempboxa}%
   \else  \@circ\@tempdimb{96}\fi\endgroup}

\def\@dot#1{\@tempdimb #1\unitlength \@circ\@tempdimb{112}}

\def\@circ#1#2{\@tempdima #1\relax \advance\@tempdima .5pt\relax
   \@tempcnta\@tempdima \@tempdima 1pt\relax
   \divide\@tempcnta\@tempdima
   \ifnum\@tempcnta > 15\relax \@tempcnta 15\relax \fi
   \ifnum \@tempcnta >\z@ \advance\@tempcnta\m@ne\fi
   \advance\@tempcnta #2\relax
   \@circlefnt \char\@tempcnta}

\def\@nnil{\@nil}
\def\@empty{}
\def\@fornoop#1\@@#2#3{}

\def\@tfor#1:=#2\do#3{\xdef\@fortmp{#2}\ifx\@fortmp\@empty \else
    \@tforloop#2\@nil\@nil\@@#1{#3}\fi}
\def\@tforloop#1#2\@@#3#4{\def#3{#1}\ifx #3\@nnil
       \let\@nextwhile=\@fornoop \else
      #4\relax\let\@nextwhile=\@tforloop\fi\@nextwhile#2\@@#3{#4}}

\def\@height{height}
\def\@depth{depth}
\def\@width{width}

\def\newbox{\alloc@4\box\chardef\insc@unt}

\def\@makebox[#1]{\leavevmode\@ifnextchar [{\@imakebox[#1]}{\@imakebox[#1][x]}}

\long\def\@imakebox[#1][#2]#3{\hbox to#1{\let\mb@l\hss
\let\mb@r\hss \expandafter\let\csname mb@#2\endcsname\relax
\mb@l #3\mb@r}}

\def\@makepicbox(#1,#2){\leavevmode\@ifnextchar
   [{\@imakepicbox(#1,#2)}{\@imakepicbox(#1,#2)[]}}

\long\def\@imakepicbox(#1,#2)[#3]#4{\vbox to#2\unitlength
   {\let\mb@b\vss \let\mb@l\hss\let\mb@r\hss
    \let\mb@t\vss
    \@tfor\@tempa :=#3\do{\expandafter\let
        \csname mb@\@tempa\endcsname\relax}%
\mb@t\hbox to #1\unitlength{\mb@l #4\mb@r}\mb@b}}

\def\newsavebox#1{\@ifdefinable#1{\newbox#1}}

\def\savebox#1{\@ifnextchar ({\@savepicbox#1}{\@ifnextchar
     [{\@savebox#1}{\sbox#1}}}

\def\sbox#1#2{\setbox#1\hbox{#2}}

\def\@savebox#1[#2]{\@ifnextchar [{\@isavebox#1[#2]}{\@isavebox#1[#2][x]}}

\long\def\@isavebox#1[#2][#3]#4{\setbox#1 \hbox{\@imakebox[#2][#3]{#4}}}

\def\@savepicbox#1(#2,#3){\@ifnextchar
   [{\@isavepicbox#1(#2,#3)}{\@isavepicbox#1(#2,#3)[]}}

\long\def\@isavepicbox#1(#2,#3)[#4]#5{\setbox#1 \hbox{\@imakepicbox
     (#2,#3)[#4]{#5}}}
\long\def\frame#1{\leavevmode
    \hbox{\hskip-\@wholewidth
     \vbox{\vskip-\@wholewidth
            \hrule \@height\@wholewidth
          \hbox{\vrule \@width\@wholewidth #1\vrule \@width\@wholewidth}\hrule
           \@height \@wholewidth\vskip -\@halfwidth}\hskip-\@wholewidth}}

\newdimen\fboxrule
\newdimen\fboxsep

\long\def\fbox#1{\leavevmode\setbox\@tempboxa\hbox{#1}\@tempdima\fboxrule
    \advance\@tempdima \fboxsep \advance\@tempdima \dp\@tempboxa
   \hbox{\lower \@tempdima\hbox
  {\vbox{\hrule \@height \fboxrule
          \hbox{\vrule \@width \fboxrule \hskip\fboxsep
          \vbox{\vskip\fboxsep \box\@tempboxa\vskip\fboxsep}\hskip
                 \fboxsep\vrule \@width \fboxrule}
                 \hrule \@height \fboxrule}}}}

\def\framebox{\@ifnextchar ({\@framepicbox}{\@ifnextchar
     [{\@framebox}{\fbox}}}

\def\@framebox[#1]{\@ifnextchar [{\@iframebox[#1]}{\@iframebox[#1][x]}}

\long\def\@iframebox[#1][#2]#3{\leavevmode
  \savebox\@tempboxa[#1][#2]{\kern\fboxsep #3\kern\fboxsep}\@tempdima\fboxrule
    \advance\@tempdima \fboxsep \advance\@tempdima \dp\@tempboxa
   \hbox{\lower \@tempdima\hbox
  {\vbox{\hrule \@height \fboxrule
          \hbox{\vrule \@width \fboxrule \hskip-\fboxrule
              \vbox{\vskip\fboxsep \box\@tempboxa\vskip\fboxsep}\hskip
                  -\fboxrule\vrule \@width \fboxrule}
                  \hrule \@height \fboxrule}}}}

\def\@framepicbox(#1,#2){\@ifnextchar
   [{\@iframepicbox(#1,#2)}{\@iframepicbox(#1,#2)[]}}

\long\def\@iframepicbox(#1,#2)[#3]#4{\frame{\@imakepicbox(#1,#2)[#3]{#4}}}

\def\parbox{\@ifnextchar [{\@iparbox}{\@iparbox[c]}}

\long\def\@iparbox[#1]#2#3{\leavevmode \@pboxswfalse
   \if #1b\vbox
     \else \if #1t\vtop
              \else \ifmmode \vcenter
                        \else \@pboxswtrue $\vcenter
                     \fi
           \fi
    \fi{\hsize #2\@parboxrestore #3}\if@pboxsw $\fi}

\catcode`\@=13

\def\Boxum#1{\underbar{#1}}
\def\GBoxum#1{\underbar{#1}}
\def\DBoxum#1{\raise-2.5pt\hbox{\framebox(8,10){\framebox(6,8){#1}}}}
\def\Kreisum#1{\kern1pt#1\kern-2.5pt\raise3pt\hbox{\circle8}\kern-5pt}
\def\GKreisum#1{\kern3pt#1\kern-4.5pt\raise3pt\hbox{\circle{11}}\kern-5pt}
\def\DKreisum#1{\kern3pt#1\kern-2.5pt\raise3pt\hbox{\circle{11}\kern-11pt
   \circle8}\kern-4pt}
\def\GDKreisum#1{\kern3pt#1\kern-4.5pt\raise3pt\hbox{\circle{14}\kern-14pt
   \circle{11}}\kern-3pt}
\def\Ferrers{\PfadDicke{.3pt}
\Pfad(0,-4),22222222\endPfad
\Pfad(2,-4),22222222\endPfad
\Pfad(4,-4),22222222\endPfad
\Pfad(6,-2),222222\endPfad
\Pfad(8,2),22\endPfad
\Pfad(0,4),11111111\endPfad
\Pfad(0,2),11111111\endPfad
\Pfad(0,0),111111\endPfad
\Pfad(0,-2),111111\endPfad
\Pfad(0,-4),1111\endPfad }
\def\Eintr#1(#2,#3){\Label\o{\raise-8pt\hbox{#1}}(#2,#3)}
\def\Zeilea#1#2#3#4{
\Eintr{#1}(1,2)
\Eintr{#2}(3,2)
\Eintr{#3}(5,2)
\Eintr{#4}(7,2)
}
\def\Zeileb#1#2#3{
\Eintr{#1}(1,0)
\Eintr{#2}(3,0)
\Eintr{#3}(5,0)
}
\def\Zeilec#1#2#3{
\Eintr{#1}(1,-2)
\Eintr{#2}(3,-2)
\Eintr{#3}(5,-2)
}
\def\Zeiled#1#2{
\Eintr{#1}(1,-4)
\Eintr{#2}(3,-4)
}
\def\Tableau#1/#2/#3/#4/{
\PfadDicke{.3pt}
\Zeilea#1
\Zeileb#2
\Zeilec#3
\Zeiled#4
\hbox{\hskip1.8cm}
}

\topmatter
\title Another involution principle-free bijective proof of Stanley's
hook-content formula
\endtitle
\author C.~Krattenthaler
\endauthor
\affil
Institut f\"ur Mathematik der Universit\"at Wien,\\
Strudlhofgasse 4, A-1090 Wien, Austria.\\
e-mail: KRATT\@Pap.Univie.Ac.At\\
WWW: \tt http://radon.mat.univie.ac.at/People/kratt
\endaffil
\address Institut f\"ur Mathematik der Universit\"at Wien,
Strudlhofgasse 4, A-1090 Wien, Austria.
\endaddress
\abstract
Another bijective proof of Stanley's hook-content formula for the generating
function for semistandard tableaux of a given shape
is given that does not involve the involution principle of Garsia and
Milne. It is the result of a merge of the modified jeu de taquin idea
from the author's previous bijective proof (``An involution
principle-free bijective proof of Stanley's hook-content formula",
Discrete Math\. Theoret\. Computer Science, to appear) 
and the Novelli-Pak-Stoyanovskii bijection 
(Discrete Math\. Theoret\. Computer Science {\bf 1} (1997), 
53--67) for the
hook formula for standard Young tableaux of a given shape. This new
algorithm can also be used as an
algorithm for the random generation of tableaux of a given shape with
bounded entries. An appropriate deformation of this algorithm
gives an algorithm for the random generation of plane partitions
inside a given box.

\endabstract
\endtopmatter
\document

\subhead 1. Introduction\endsubhead
There are two recent major contributions to bijective combinatorics in
regard of bijective proofs of hook formulas for various types of
tableaux: First, Novelli, Pak and Stoyanovskii \cite{\NoPSAA}
(announced in \cite{\PaStAA}) came up with a new bijective proof of
the celebrated Frame--Robinson--Thrall hook formula \cite{\FrRTAA} for the
number of standard Young tableaux of a given shape, which is probably
{\it the} bijective proof of the hook formula that everybody was looking for.
This bijection can also be used as an algorithm for
the random generation of standard Young tableaux of a given
shape. Second, in \cite{\KratAY}
the author of this paper gave the first bijective proof of
Stanley's hook-content formula \cite{\StanAA, Theorem~15.3} for the
generating function for semistandard tableaux of a given shape with
bounded entries which does not rely on the involution principle of Garsia and
Milne \cite{\GaMiAB}.

The purpose of this paper is to merge the ideas from the
aforementioned bijections and thus obtain another bijective proof of
Stanley's hook-content formula. Unlike the previous bijection, the
new bijection can also be used as an algorithm for
the random generation of semistandard tableaux of a given
shape with bounded entries. A simple deformation of this algorithm in
the rectangular shape case yields, as well, an algorithm for
the random generation of plane partitions inside a given box.

In order to
be able to state the formula, we have to recall some basic definitions.
A {\it partition\/} is a sequence
$\la=(\la_1,\la_2,\dots,\la_r)$ with $\la_1\ge \la_2\ge\dots\ge\la_r>0$,
for some $r$. The {\it Ferrers
diagram\/} of $\lambda$
is an array of cells with $r$ left-justified rows and $\lambda_i$
cells in row $i$. Figure~1.a shows the Ferrers
diagram corresponding to $(4,3,3,2)$.
The {\it conjugate of\/} $\lambda$ is the partition
$(\lambda^\prime_1, \dots, \lambda^\prime_{\lambda_1})$ where
$\lambda_j'$ is the length of the $j$-th column in the Ferrers diagram of
$\lambda$. We label the cell in the $i$-th row and $j$-th column of
(the Ferrers diagram of) $\la$ by the pair $(i,j)$. Also, if we write
$\rh\in\la$ we mean `$\rh$ is a cell of $\la$'. The {\it hook length\/}
$h_\rh$ of a cell $\rh=(i,j)$ of $\la$ is $(\la_i-j)+(\la'_j-i)+1$, \
the number of cells in the
{\it hook\/} of $\rh$, which is the set of cells that are either in
the same row as $\rh$ and to the right of $\rh$, or in the same column
as $\rh$ and below $\rh$, $\rh$ included. The hook of a cell
$\rh=(i,j)$ consists of three parts, the cell itself, the {\it arm\/},
which is the set of cells to the right of $\rh$ in the same row as
$\rh$, and the {\it leg\/}, which is the set of cells below
$\rh$ in the same column as $\rh$. The number of cells in the arm of
$\rh$ (which is $\la_i-j$) is denoted by $a_\rh$, and
the number of cells in the leg of
$\rh$ (which is $\la_j'-i$) is denoted by $l_\rh$.
The {\it content\/} $c_\rh$
of a cell $\rh=(i,j)$ of $\la$ is $j-i$.
\vskip10pt
\vbox{
$$
\gather
\smatrix \format \sa\c\quad \s\c\quad \s\c\quad \s\c\quad \se\\
\hlinefor9\\
&\vphantom{f}&& && && &\\
\hlinefor9\\
&\vphantom{f}&& && &\\
\hlinefor7\\
&\vphantom{f}&& && &\\
\hlinefor7\\
&\vphantom{f} && & \\
\hlinefor5
\endsmatrix
\hskip1cm
\hskip1cm
\smatrix \format\sa\c\s\c\s\c\s\c\se\\
\hlinefor9\\
&\hbox to10pt{\hss 1\hss}&&\hbox to10pt{\hss 1\hss}&&\hbox to10pt{\hss
2\hss}&&\hbox to10pt{\hss 5\hss}&\\
\hlinefor9\\
&2&&4&&6&\\
\hlinefor7\\
&4&&5&&7&\\
\hlinefor7\\
&5&&6&\\
\hlinefor5
\endsmatrix
\hskip1cm
\hskip1cm
\smatrix \format\sa\c\s\c\s\c\s\c\se\\
\hlinefor9\\
&\hbox to10pt{\hss 1\hss}&&\hbox to10pt{\hss 5\hss}&&\hbox to10pt{\hss
9\hss}&&\hbox to10pt{\hss 12\hss}&\\
\hlinefor9\\
&2&&6&&10&\\
\hlinefor7\\
&3&&7&&11&\\
\hlinefor7\\
&4&&8&\\
\hlinefor5
\endsmatrix\\
\smallmatrix\format\l\\
 \text{\eightpoint a. Ferrers diagram\quad \quad \quad \quad \quad
\quad
b. semistandard tableau
\quad \quad \quad \quad \quad
c. order of the cells}\\%
\endsmallmatrix
\endgather
$$
\centerline{\eightpoint Figure 1}
}
\vskip10pt

Given a partition $\la=(\la_1,\la_2,\dots,\la_r)$, a {\it semistandard
tableau} ({\it tableau} for short) {\it of shape\/} $\la$ is a filling
$T$ of the cells of $\la$
with integers such that the entries along rows are weakly increasing
and such that the entries along columns are
stricly increasing. Figure~1.b displays a (semistandard) tableau
of shape $(4,3,3,2)$.
We write $T_\rh$ for the entry in cell $\rh$ of $T$. We
call the sum of all the entries of a tableau $T$ the
{\it norm\/} of $T$, and denote it by $n(T)$.

Now we are in the position to state Stanley's hook-content formula
\cite{\StanAA, Theorem~15.3}.

\proclaim{Theorem 1 (Stanley)}Let $\la=(\la_1,\la_2,\dots,\la_r)$ be a
partition and $b$ be an integer $\ge r$. The generating function
$\sum _{} ^{} q^{n(T)}$, where the sum is over all semistandard
tableaux $T$ of shape $\la$ with entries between $1$ and $b$, is given by
$$q^{\sum _{i=1} ^{r}i\la_i}\prod _{\rh\in\la} ^{}\frac
{1-q^{b+c_\rh}} {1-q^{h_\rh}}.\tag1.1$$
\endproclaim

We refer the reader to the Introduction of \cite{\KratAY} for a detailed
exposition on proofs of this formula. The bijective proofs which are
known for (1.1) are the rather involved proof by Remmel and Whitney
\cite{\ReWhAA},
which makes use of the involution principle, and the proof by the
author \cite{\KratAY}, which makes use of the Hillman-Grassl algorithm
\cite{\HiGrAA}
and a modified jeu de taquin (the jeu de taquin being due to
Sch\"utzenberger \cite{\SchuAA}).
The latter bijection is actually a bijection for
a rewriting of (1.1), where the ``content product" is brought to the
other side,
$$\Big(\underset \text{with }1\le \text{entries}\le b\to {\sum _{T
\text{ a tableau of
shape $\la$}}} ^{}q^{n(T)}\Big)\cdot\prod _{\rh\in\la} ^{}\frac {1}
{1-q^{b+c_\rh}}
= q^{\sum _{i=1} ^{r}i\la_i}\prod _{\rh\in\la} ^{}\frac {1}
{1-q^{h_\rh}}.
\tag1.2$$
(It is shown in Section~3 of \cite{\KratAY} that the same ideas also
yield a bijective proof for (1.1), directly.)
On the other hand, Remmel and Whitney's bijection is one for a
rewriting of (1.1), where the ``hook product" is brought to the other
side, and where each factor of the form $(1-q^\al)$ is divided by
$(1-q)$, so, basically for
$$\Big(\underset \text{with }1\le \text{entries}\le b\to {\sum _{T
\text{ a tableau of
shape $\la$}}} ^{}q^{n(T)}\Big)\cdot\prod _{\rh\in\la} ^{}
(q^{-a_\rh}+q^{-a_\rh+1}+\dots+q^{l_\rh})
= \prod _{\rh\in\la} ^{}
(q^{1-c_\rh}+q^{2-c_\rh}+\dots+q^{b}).
\tag1.3$$
In \cite{\KratAY, Sec.~3} a bijection for this rewriting was
constructed which is simpler than the bijection by Remmel and Whitney,
using the basic algorithm from the bijection for (1.2) from the same
paper. However, we did have to make use of the involution principle.

In the next section, we describe a new bijection for (1.3), which does
not rely on the involution principle. It is a merge of the modified
jeu de taquin idea from \cite{\KratAY} and the
Novelli-Pak-Stoyanovskii bijection \cite{\NoPSAA} for the
hook formula for standard Young tableaux of a given shape.
A few properties of our algorithm, which are crucial in showing that
it is indeed the desired bijection, are stated and proved separately
in Section~3.
Finally, in Section~4, we argue that our bijection can be used as an
algorithm for the random generation of tableaux of a given shape with
bounded entries, and we also work out the appropriate deformation
explicitly
which gives an algorithm for the random generation of plane partitions
inside a given box.

\subhead 2. Bijective proof of Stanley's hook-content
formula\endsubhead
In this section we present a bijective proof of Theorem~1, by
constructing a bijection for the equivalent identity (1.3).
In all what follows, we assume that we are given a fixed partition
$\la=(\la_1,\la_2,\dots,\la_r)$
and a fixed positive integer $b$ which is at least as large as
$r$, the number of rows of the Ferrers diagram of $\la$.

In order to formulate (1.3) combinatorially, it is convenient
to introduce a few terms. We
call a filling of the cells of $\la$ with
integers such that the entry in cell $\rh$ is at least $-a_\rh$ and at
most $l_\rh$, for any cell $\rh$ of $\la$, a {\it hook
tabloid of shape\/} $\la$. Furthermore, we
call a filling of the cells of $\la$ with
integers such that the entry in cell $\rh$ is at least $1-c_\rh$ and at
most $b$, for any cell $\rh$ of $\la$, a {\it content
tabloid of shape\/} $\la$. We use the terms that we introduced for
tableaux also for (either kind of) tabloids. In particular,
the sum of the entries of some tabloid $T$ is again called the {\it
norm\/} of $T$ and denoted by $n(T)$.

In terms of these notions, identity (1.3) can be rewritten as
$$\sum _{(\hT,\hH)} ^{}q^{n({\hT})}q^{n({\hH })}=
\sum _{\hC} ^{}q^{n({\hC})},$$
where the left-hand side sum is over all
pairs $({\hT},{\hH })$, with ${\hT}$ varying over all
tableaux of shape $\la$ with entries between $1$ and $b$ and
${\hH }$ varying over all hook tabloids of shape $\la$, and where
the right-hand side sum 
is over all content tabloids $\hC$ of shape $\la$. So the task is to set
up a bijection that maps a right-hand side content tabloid ${\hC}$ to a
left-hand side pair $({\hT},{\hH })$, such that
$$n({\hT})+n({\hH })=n({\hC}).\tag2.1$$

We claim that Algorithm~HC, which is described below, performs this task.
For the formulation of the algorithm, it is convenient to fix, once and for
all, a particular total order of the cells of $\la$. In this order,
the cells in each column are ordered from top to bottom, and cells
in column 1 come before cells in column 2, etc. Figure~1.c
illustrates this particular order of the cells for $\la=(4,3,3,2)$.

The building blocks of Algorithm~HC are jeu de taquin-like moves,
which are explained in Definition~JT below.

\definition{Definition JT} Consider a filling of shape $\la$ with
integers, and a cell $\om$ filled with entry $s$. Suppose that the
entries which are in the region consisting of all the successors
of $\om$, in the fixed total order of the cells,
are weakly increasing along
rows and strictly increasing along columns.

We say that we {\it move $s$ to the right/bottom by modified jeu de
taquin}, when we perform the algorithm:
\smallskip
Call $s$ the special entry.
\smallskip
(JT) Compare the special entry $s$ with its right neighbour,
$x$ say (if there is no right neighbour, then, by convention, 
we set $x=\infty$), and its bottom neighbour, $y$ say 
(if there is no bottom neighbour, then, by convention, we set $y=\infty$), 
see (2.2). If $s\le x$ and $s<y$, then stop.

If not, i.e\. if the special entry violates weak increase along
rows or strict increase along columns, then we have the following
situation,
$$\smatrix \format \sa\c\s\c\se\\
\hlinefor5\\
&s&&x&\\
\hlinefor5\\
&y&\\
\hlinefor3
\endsmatrix\ ,\tag2.2
$$
where at least one of $s>x$ and $s\ge y$ holds. If $x+1<y$ then do the move
$$\smatrix \format \sa\c\s\c\se\\
\hlinefor5\\
&x+1&&s&\\
\hlinefor5\\
&y&\\
\hlinefor3
\endsmatrix\ .\tag2.3
$$
If $x\ge y-1$ then do the move
$$\smatrix \format \sa\c\s\c\se\\
\hlinefor5\\
&y-1&&x&\\
\hlinefor5\\
&s&\\
\hlinefor3
\endsmatrix\ .\tag2.4
$$
The ``new" special entry in either case is $s$ again.
Repeat (JT). (Note that always after either type of
move the only possible violations of increase along rows or strict
increase along columns in the region consisting of cell $\om$ and all
its successors in the fixed total order of the cells
involve the special entry and the entry to
the right or/and below.)\quad \quad \qed

\enddefinition

Clearly, aside from adding/subtracting $1$ to/from
the special entry, what happens from (2.2) to (2.3), respectively
(2.4), is a {\it jeu de taquin forward move\/} (cf\. \cite{\SchuAA,
Sec.~2; \SagaAL, pp.~120/169}).
Because of that, we call the way
the special entry takes during the performance of (JT) the {\it
forward\/} ({\it jeu de taquin}) {\it path of\/} $\om$, the cell
where $s$ starts.

\example{Example JT} Figure~2 shows an example for the algorithm
defined in Definition~JT, with $\la=(6,5,5,4)$ and $\om=(2,2)$. The array
in the left half of Figure~2 displays a filling of $(6,5,5,4)$ which
does have the property that entries in the region consisting of
all the successors of $(2,2)$ in the fixed total order of cells
are weakly increasing along rows
and strictly increasing along columns. The
entry in cell $(2,2)$ is circled. The right half of Figure~2 displays
the array which is obtained after application of the algorithm. The path
which the special entry takes is marked by arrows, and the special
entry is circled
in the cell where it finally stops.
\vskip10pt
\vbox{
$$
\smatrix \format \sa\r\s\r\s\r\s\r\s\r\s\r\se\\
\hlinefor{13}\\
&4&&1&&2&&2&&3&&4&\\
\hlinefor{13}\\
&4&&\Kreisum6&&3&&4&&4&\\
\hlinefor{11}\\
&8&&5&&5&&5&&6&\\
\hlinefor{11}\\
&6&&6&&7&&7&\\
\hlinefor9
\endsmatrix
\quad \overset\text {(JT)}\to\longrightarrow\quad 
\smatrix \format \sa\r\s\r\s\r\s\r\s\r\s\r\se\\
\hlinefor{13}\\
&4&&1&&2&&2&&3&&4&\\
\hlinefor{13}\\
&4&&\hbox to5pt{$4\kern-4pt\rightarrow$\hss}&&4&&4&&4&\\
\hlinefor{11}\\
&8&&5&&
\vbox to10pt{\hsize0pt \vss \noindent\hbox to0pt{$\uparrow$\hss}\newline\vphantom{f}}
\hbox to5pt{$6\kern-4pt\rightarrow$\hss}&&\Kreisum6&&6&\\
\hlinefor{11}\\
&6&&6&&7&&7&\\
\hlinefor9
\endsmatrix
$$
\centerline{\eightpoint Figure 2}
}
\vskip10pt
\line{\hfill$\square$}
\endexample

Now, the idea of Algorithm~HC is to start with a content tabloid $\hC$, 
and convert it, step by step, into a tableau $\hT$. This
conversion uses the above modifed jeu de taquin. The moves of the entries
during this modified jeu de taquin are recorded, in an appropriate way, by a
hook tabloid $\hH $.

\bigskip
\noindent
{\bf Algorithm HC}. The input for the algorithm is a content tabloid
${\hC}$ of shape $\la$.

\smallskip
(HC0) Set $(T,H):=({\hC},\bold 0)$, where $\bold 0$ denotes the hook tabloid
of shape $\la$ with $0$ in each cell.
Call the last cell of $\la$ in the fixed total order of the
cells (i.e., it is the bottommost cell in the last column)
the {\it distinguished cell}. Continue with (HC1).

\smallskip
(HC1) Denote the distinguished cell by $\om$.
Move the entry in the distinguished cell
to the right/bottom by modified jeu de
taquin according to Definition~JT. 
Say that it stops in cell $\om'$. Continue with (HC2).

\smallskip
(HC2) Let $T$ be the array just
obtained. Let $\om=(i,j)$ and $\om'=(i',j')$. 
The ``new" $H$ is obtained from the ``old"
$H$ by moving all the entries of the old $H$
in $\om$'s column between row $i+1$ and $i'$
up by one row, thereby also adding 1 to each such entry, and by
setting the entry in cell $(i',j)$ equal to $(j-j')$, the negative of
the number of columns we moved to the right during (HC1). 
In symbols: For $k=i,i+1,\dots,i'-1$, the
entry $H_{(k,j)}$ is replaced by $H_{(k+1,j)}+1$, and $H_{(i',j)}$ is
replaced by $(j-j')$. All other entries are left unchanged.

If the distinguished cell is the first cell of $\la$ 
(i.e., the topmost cell in the first
column), then stop. The output of the algorithm is $(T,H)$.

Otherwise, we choose as the
``new" distinguished cell the predecessor of the ``old" distinguished
cell in the fixed total
order of the cells. Continue with (HC1).\quad \quad \qed

\bigskip
\example{Example HC}A complete example for Algorithm~HC can be found in
the appendix. There we choose $b=7$ and map the content tabloid
of shape $(4,3,3,2)$ on the left of Figure~3 to the pair on
the right of Figure~3, consisting of a tableau
of shape $(4,3,3,2)$ with entries between 1 and 7 
and a hook tabloid of shape $(4,3,3,2)$, such that the
weight property (2.1) is satisfied. In fact, the norm of the content tabloid
on the left of Figure~3 is $46$, while the norm of the
tableau on the right is $48$ and the
norm of the hook tabloid is $-2$.
\vskip10pt
\vbox{
$$
\matrix \format \c&\quad \c\quad &\c\\
\smatrix \format \sa\r\s\r\s\r\s\r\se\\
\hlinefor9\\
&\hphantom{-{}}7&&\hphantom{-{}}3&&\hphantom{-{}}5&&-2&\\
\hlinefor9\\
&7&&3&&2&\\
\hlinefor7\\
&5&&4&&2&\\
\hlinefor7\\
&4&&6&\\
\hlinefor5
\endsmatrix
&\longleftrightarrow&
\(
\smatrix \format \sa\r\s\r\s\r\s\r\se\\
\hlinefor9\\
&\hphantom{-{}}1&&\hphantom{-{}}1&&\hphantom{-{}}2&&\hphantom{-{}}5&\\
\hlinefor9\\
&2&&4&&6&\\
\hlinefor7\\
&4&&5&&7&\\
\hlinefor7\\
&5&&6&\\
\hlinefor5
\endsmatrix
,
\smatrix \format \sa\r\s\r\s\r\s\r\se\\
\hlinefor9\\
&3&&0&&-1&&\hphantom{-{}}0&\\
\hlinefor9\\
&-1&&-1&&1&\\
\hlinefor7\\
&-2&&-1&&0&\\
\hlinefor7\\
&0&&0&\\
\hlinefor5
\endsmatrix
\)\\
n(.)=46& &n(.)=48\ ,\ n(.)=-2
\endmatrix
$$
\centerline{\eightpoint Figure 3}
}
\vskip10pt
The appendix has to be read in the following way. What the
left columns show is the pair $(T,H)$ that is obtained after each loop
(HC1)-(HC2). Together with the pair $(T,H)$ a picture of the shape
$(4,3,3,2)$ is displayed, containing certain paths.
This will be important for understanding the inverse algorithm
(Algorithm~HC* below) but can be ignored for
the moment.
At each stage, the entry in the distinguished cell is circled.
Then each modified jeu de taquin move during the execution of (JT) in 
(HC1) is displayed in the
right columns. The special entry is always
underlined. Once the special entry stops,
the filling of $(4,3,3,2)$ that is obtained is displayed again
in the left column of the row below, with the special entry
doubly circled. In the filling, the portion consisting of all the
successors of the distinguished cell in the fixed total order of the
cells (this is the portion which is ``already ordered") is separated
from the rest by a thick line. The hook tabloid is
changed according to the rules described in (HC2) and is displayed
next to the filling.\quad \quad \qed
\endexample

We have to prove that Algorithm~HC is well-defined, that the weight
property (2.1) is satisfied, and that it is a bijection between
``right-hand side objects" and ``left-hand side objects" of (1.3).

For Algorithm~HC being well-defined, we have to show that the
output is always a pair $(\hT,\hH )$, where $\hT$ is a
tableau of shape $\la$ with entries between $1$ and $b$, and where
${\hH }$ is a hook tabloid of shape $\la$.

That $\hT$ is indeed a tableau follows from the observation that
after each loop (HC1)-(HC2) the portion of the array $T$, obtained
after that loop, consisting of cell $\om$ and all cells which are after
$\om$ in the fixed total order of the cells, is a {\it skew
tableau\/}, i.e., the entries are weakly increasing along rows and
strictly increasing along columns in that portion.
This observation is based on the easily verified fact
that always after either type of
move, (2.2)-(2.3) or (2.2)-(2.4),
the only possible violations of increase along rows or strict
increase along columns in that portion
involve the special entry and the entry to
the right or/and below, as we already mentioned.

Furthermore, the entries of $\hT$ are at most $b$. For, in the
beginning we start with $T=\hC$, which is an array of integers all
of which are at most $b$. Then, in order to finally obtain $\hT$,
we perform moves (2.2)-(2.3) and (2.2)-(2.4). We claim that at any
stage, the entries will be at most $b$. The only possible
problem could arise from the move (2.2)-(2.3), where $x$ is increased
by 1. But, by induction, we may assume that all entries are $\le b$
before such a move. In particular, we have $s\le b$ and $y\le b$.
Now, the move (2.2)-(2.3) is applied only if $x+1<y$, or if $y$ is
not there and, hence, $s>x$. In either case it follows that $x+1\le
b$.

We claim that the entries of $\hT$ are also at least $1$. Since we
already know that $\hT$ is a tableau, it suffices to show that the
entry, $e$ say, in cell $(1,1)$ is at least $1$. Now, the type of
moves (2.2)-(2.3) and (2.2)-(2.4) imply that the origin of the entry
$e$ must be some entry $e'$ in cell $\rh$ of $\hC$, which was
moved to the left and up according to (2.3) and (2.4), thereby
adding/subtracting 1, in a way such that $e=e'+c_\rh$. Since $\hC$ is a
content tabloid, we have $e'\ge 1-c_\rh$. Thus, it follows that
$e=e'+c_\rh\ge 1$, as desired.

This proves that $\hT$ is indeed a
tableau of shape $\la$ with entries between $1$ and $b$.

\smallskip
Next, $\hH $ is a hook tabloid, because the way the new array $H$
is formed in (HC2) implies directly that it is a hook tabloid at any
stage, and so must be $\hH $ at the end.

\smallskip
Similarly, the weight property (2.1) holds because the algorithm is
constructed in a way that it holds at any stage. To be precise,
whenever we performed (HC2), we have
$$n({T})+n({H})=n({\hC}).\tag2.5$$
For, equation (2.5) trivially holds at the very beginning (where we
set $T=\hC$ and $H=\bold 0$). Then, during each loop (HC1)-(HC2),
the special entry $s$ is moved from cell $\om=(i,j)$ to cell
$\om'=(i',j')$. In each loop, while applying (JT) (possibly
repeatedly), we add 1 to $j'-j$
entries of $T$, and we subtract 1 from $i'-i$ entries of $T$.
On the other hand, in
(HC2) we add 1 to $i'-i$ entries of $H$, and we introduce a new
entry $j-j'$ into $H$. Obviously, these alterations add up to zero.
This establishes (2.5), and thus, by spezializing (2.5) to the final
output $(\hT,\hH )$ of Algorithm~HC, equation (2.1) as well.

\smallskip
What remains is to demonstrate that our algorithm is actually a
{\it bijection\/} between right-hand side and left-hand side objects
of (1.3). This
will be accomplished by constructing an algorithm, Algorithm~HC* below,
that will turn out to be the inverse of Algorithm~HC. That
Algorithm~HC* is well-defined will be argued after the statement of
the algorithm, with two critical properties being established
separately in Lemmas~HC* and OC* in Section~3, respectively. 
That Algorithm~HC and Algorithm~HC* are
inverses of each other will follow from the construction of the
algorithms and from Lemma~HC in Section~3.

\smallskip
As the reader may anticipate, the building blocks of Algorithm~HC*
are again jeu de taquin-like moves, which are, of course, inverse
to those in Definition~JT. We explain them in Definition~JT* below.

\definition{Definition JT*} Consider a filling of shape $\la$ with
integers, a cell $\om$, and a cell $\om'$, filled with entry $s$,
which is located weakly to the right and below of $\om$. Suppose that the
entries which are in the region consisting of cell $\om$ and all
its successors in the fixed total order of the cells
are weakly increasing along rows and strictly increasing along columns.

We say that we {\it move $s$ to the left/top until $\om$ by modified jeu de
taquin}, when we perform the algorithm:
\smallskip
Call $s$ the special entry.
\smallskip
(JT*) If the special entry $s$ is located in $\om$, then stop.

If not, then we have the following
situation,
$$\smatrix \format \sa\c\s\c\se\\
\omit&\omit&\hlinefor3\\
\omit&&&y&\\
\hlinefor5\\
&x&&s&\\
\hlinefor5
\endsmatrix\ ,\tag2.6
$$
(If $s$ should be in the same column as $\om$, by convention
we set $x=\infty$. If $y$ is
actually not there, by convention we set $y=\infty$.) 
If $x-1>y$ then do the move
$$\smatrix \format \sa\c\s\c\se\\
\omit&\omit&\hlinefor3\\
\omit&&&y&\\
\hlinefor5\\
&s&&x-1&\\
\hlinefor5
\endsmatrix\ .\tag2.7
$$
If $x\le y+1$ then do the move
$$\smatrix \format \sa\c\s\c\se\\
\omit&\omit&\hlinefor3\\
\omit&&&s&\\
\hlinefor5\\
&x&&y+1&\\
\hlinefor5
\endsmatrix\ .\tag2.8
$$
The ``new" special entry in either case is $s$ again.
Repeat (JT*).\quad \quad \qed

\enddefinition

Again, it should be noticed that what happens from (2.6) to (2.7),
respectively from (2.6) to (2.8), are {\it jeu de
taquin backward moves\/} (cf\. \cite{\SchuAA, Sec.~2; \SagaAL,
pp.~120/169}),
which reverse the forward moves (2.2)-(2.3)
and (2.2)-(2.4), respectively.
Because of that, we call
the way that the special entry $s$ takes during the execution of (JT*) the
{\it backward\/} ({\it jeu de taquin}) {\it path of $\om'$}
(the cell where $s$ starts) {\it until $\om$}.

\example{Example JT*} Figure~4 shows an example for the algorithm
defined in Definition~JT*, with $\la=(6,5,5,4)$, $\om=(2,2)$ and
$\om'=(3,4)$. The array
in the left half of Figure~4 displays a filling of $(6,5,5,4)$ which
does have the property that entries in the region consisting of
cell $(2,2)$ and all its
successors in the fixed total order of cells are weakly increasing along rows
and strictly increasing along columns. The 
entry in cell $(2,2)$ is doubly circled, the entry in cell $(3,4)$ is
circled. 
The right half of Figure~4 displays
the array which is obtained after application of the algorithm. The path
which the special entry takes is marked by arrows, and the special
entry is circled
in the cell where it finally stops.
\vskip10pt
\vbox{
$$
\smatrix \format \sa\r\s\r\s\r\s\r\s\r\s\r\se\\
\hlinefor{13}\\
&4&&1&&2&&2&&3&&4&\\
\hlinefor{13}\\
&4&&\DKreisum4&&4&&4&&4&\\
\hlinefor{11}\\
&8&&5&&6&&\Kreisum6&&6&\\
\hlinefor{11}\\
&6&&6&&7&&7&\\
\hlinefor9
\endsmatrix
\quad \overset\text {(JT*)}\to\longrightarrow\quad 
\smatrix \format \sa\r\s\r\s\r\s\r\s\r\s\r\se\\
\hlinefor{13}\\
&4&&1&&2&&2&&3&&4&\\
\hlinefor{13}\\
&4&&\Kreisum6&&\hbox to5pt{\hss$\leftarrow\kern-4pt3$}&&4&&4&\\
\hlinefor{11}\\
&8&&5&&\vbox to10pt{\hsize5.27pt \vss \noindent$\uparrow$\newline$5$\vphantom{f}}&&
\hbox to5pt{\hss$\leftarrow\kern-4pt5$}&&6&\\
\hlinefor{11}\\
&6&&6&&7&&7&\\
\hlinefor9
\endsmatrix
$$
\centerline{\eightpoint Figure 4}
}
\vskip10pt
\line{\hfill$\square$}
\endexample
The reader should observe that there are actually {\it two} possible outcomes
in Definition~JT*: We may actually reach $\om$ (then the algorithm
stops there), or we may miss $\om$. As it turns out, in the
situations in which the above algorithm is applied in
Algorithm~HC* below we do always reach $\om$. This is a nontrivial fact that
is proved in Lemma~HC*.

\medskip
In order to state Algorithm~HC* we need one more definition.
In Algorithm~HC* we need to compare
different backward (jeu de taquin) paths that end in the same cell $\om$ (and
start in cells weakly to the right and below of $\om$). For that
purpose we impose a total order on (certain) backward paths as follows.

\definition{Definition P} 
Given two backward
paths $P$ and $Q$ ending in $\om$, which have the property that they
do not cross each other (i.e., they may touch each other, but never
change sides), we say that $P$ comes before $Q$ if one of the
following three conditions is satisfied:
\roster
\item"(a)" $P$ is neither contained in $Q$ nor contains $Q$ and is ``to
the right and above of $Q$", as exemplified in Figure~5, Case (a).
\item"(b)" $P$ is contained in $Q$, and $Q$ enters $P$ from below, 
as exemplified in Figure~5, Case (b).
\item"(c)" $P$ contains $Q$ and enters $Q$ from the right, 
as exemplified in Figure~5, Case (c).
\endroster

\vskip10pt
\vbox{
$$
\Einheit.4cm
\hbox to1pt{\hss}
\raise1pt\hbox{
\DickPunkt(6,-4)
\Label\ro{P}(4,-3)
\Pfad(0,-2),22\endPfad
\Pfad(0,-2),11\endPfad
\Pfad(2,-3),2\endPfad
\Pfad(2,-3),111\endPfad
\Pfad(5,-4),2\endPfad
\Pfad(5,-4),1\endPfad
}
\hbox to-1pt{\hss}
\Pfad(0,-2),22\endPfad
\Pfad(0,-2),11\endPfad
\Pfad(2,-5),222\endPfad
\Pfad(2,-5),1\endPfad
\Pfad(3,-6),2\endPfad
\DickPunkt(0,0)
\DickPunkt(3,-6)
\Label\o{\om}(0,0)
\Label\lu{Q}(2,-4)
\hskip3.5cm
\hbox to0pt{\hss}
\raise0pt\hbox{
\DickPunkt(6,-4)
\Label\ro{Q}(4,-4)
\Pfad(0,-2),22\endPfad
\Pfad(0,-2),11\endPfad
\Pfad(2,-4),22\endPfad
\Pfad(2,-4),1111\endPfad
}
\PfadDicke{2pt}
\hbox to-0pt{\hss}
\Pfad(0,-2),22\endPfad
\Pfad(0,-2),11\endPfad
\Pfad(2,-3),2\endPfad
\DickPunkt(0,0)
\DickPunkt(2,-3)
\Label\o{\om}(0,0)
\Label\u{P}(1,-2)
\hskip3.5cm
\hbox to0pt{\hss}
\raise0pt\hbox{
\PfadDicke{1pt}
\DickPunkt(6,-4)
\Label\ro{P}(4,-3)
\Pfad(0,-2),22\endPfad
\Pfad(0,-2),11\endPfad
\Pfad(2,-3),2\endPfad
\Pfad(2,-3),111\endPfad
\Pfad(5,-4),2\endPfad
\Pfad(5,-4),1\endPfad
}
\PfadDicke{2pt}
\hbox to-0pt{\hss}
\Pfad(0,-2),22\endPfad
\Pfad(0,-2),11\endPfad
\Pfad(2,-3),2\endPfad
\DickPunkt(0,0)
\DickPunkt(2,-3)
\Label\o{\om}(0,0)
\Label\u{Q}(1,-2)
\hskip2.4cm
$$
\centerline{\eightpoint
Case (a) \hskip2cm Case (b) \hskip2cm Case (c)}
\vskip3pt
\centerline{\eightpoint Figure 5}
}
\vskip10pt
\line{\hfill$\square$}
\enddefinition

Now we are in the position to formulate Algorithm~HC*.
\bigskip
\noindent
{\bf Algorithm HC*}. The input for the algorithm is a pair $(\hT,\hH )$,
where ${\hT}$ is a
tableau of shape $\la$ with entries between $1$ and $b$, and where
${\hH }$ is a hook tabloid of shape $\la$.

\smallskip
(HC*0) Set $(T,H):=(\hT,\hH )$.
Call the first cell of $\la$ in the fixed total order of the
cells (i.e., it is the topmost cell in the first column) the
{\it distinguished cell}. Continue with (HC*1).

\smallskip
(HC*1) Let the distinguished cell be $\om=(i,j)$. For each
{\it nonpositive} entry, $e$ say, of $H$ in cell $(k,j)$, $k\ge i$
(i.e., $e$ is located in the same column as $\om$ and in a row weakly
below $\om$), mark the cell $(k,j-e)$ as a {\it candidate cell\/}.

For each candidate cell determine its backward path until $\om$
according to Definition~JT*.
(In Lemma~HC* below we prove that any such backward path does indeed terminate 
in $\om$.)
Let $\om'=(i',j')$ be the candidate cell whose backward path is the first
in the total order of backward paths defined in Definition~P.
(Note that this makes indeed sense, since, clearly, two such
backward paths can never cross each other.)

Move the entry in $\om'$ to the left/top until $\om$ by modified jeu de
taquin according to Definition~JT*.
Continue with (HC*2).

\smallskip
(HC*2) Let $T$ be the array just obtained.
The ``new" $H$ is obtained from the ``old"
$H$ by moving all the entries of the old $H$
which are in the same column as $\om$, i.e., in column $j$, 
and between row $i$ and $i'-1$
down by one row, thereby also subtracting 1 from each such entry, and by
setting the entry in cell $(i,j)$ equal to $0$.
In symbols: For $k=i+1,i+2,\dots,i'$, the
entry $H_{(k,j)}$ is replaced by $H_{(k-1,j)}-1$, and $H_{(i,j)}$ is
replaced by $0$. All other entries are left unchanged.

If the distinguished cell 
is the last cell of $\la$ (i.e., the bottommost cell in the
last column), then stop. The output of the algorithm is $T$.

Otherwise, we choose as the
``new" distinguished cell the successor of the ``old" distinguished cell 
in the fixed total order of cells.
Continue with (HC*1).\quad \quad \qed

\bigskip
\example{Example HC*}A complete example for Algorithm~HC* can be found in
the appendix. There we choose $b=7$ and map the pair on
the right of Figure~3, consisting of a tableau
of shape $(4,3,3,2)$ with entries between 1 and 7 and a hook tabloid of shape
$(4,3,3,2)$, to the content tabloid
of shape $(4,3,3,2)$ on the left of Figure~3, such that the
weight property (2.1) is satisfied. It is simply the inverse of the example
for Algorithm~HC given in Example~HC. Therefore, here the appendix has to
be read in the reverse direction,
and in the following way. What the
left columns show is the pair $(T,H)$ that is obtained after each loop
(HC*1)-(HC*2), together with a copy of the Ferrers board 
$(4,3,3,2)$ in which the candidate cells are marked by bold dots, the
distinguished cell is marked by a circle, and which contains the backward paths of all
the candidate cells. The entry in the candidate cell with first
backward path is doubly circled in the tableau.
Then each
modified jeu de taquin move during the execution of (JT*) in (HC*1)
is displayed in the
right columns in the row above. The special entry is always
underlined. Once the special entry reaches the distinguished cell,
the filling of $(4,3,3,2)$ that is obtained is displayed again
in the left column, with the special entry
circled. In the filling, the portion consisting of all the
successors of the distinguished cell in the fixed total order of the
cells (this is the portion which is ``still ordered") is separated
from the rest by a thick line. The hook tabloid is
changed according to the rules described in (HC*2) and is displayed
next to the filling.\quad \quad \qed
\endexample

We have to prove that Algorithm~HC* is well-defined, that the weight
property (2.1) is satisfied, and that it is exactly the inverse of
Algorithm~HC.

In order to show that Algorithm~HC* is always well-defined, we have to
show that whenever we are in (HC*1) the 
backward path of {\it any} candidate cell reaches $\om$,
that when performing the changes to $H$ in 
(HC*2) we always obtain a hook tabloid, and, finally,
that the output is a content tabloid. 

\smallskip
The first of these three assertions is proved
in Lemma~HC* in Section~3, the third is proved in Lemma~OC* in
Section~3.

\smallskip
To see that $H$ is always a hook tabloid while performing
Algorithm~HC*, we argue by induction on the number of loops
(HC*1)-(HC*2). We have to prove that for the entry $H_\rh$ in cell
$\rh$ we have $-a_\rh\le H_\rh\le l_\rh$.
In view of the definition of transformations on $H$ in (HC*2), 
the only problem that can occur is that the
lower bound, $-a_\rh$, is violated for some entry $H_\rh$. 
We will therefore concentrate on that problem.

Trivially, $H$ is a hook tabloid at the very beginning
because in (HC*0) we set $H=\hH$. Now assume that we are at some
stage during the execution of Algorithm~HC*, that we obtained $(T,H)$
so far, where $H$ is a hook tabloid, and that we have to perform
(HC*1) with distinguished cell $\om=(i,j)$ next. 
There we determine the candidate cells, and the
candidate cell, $\om'=(i',j')$ say, with first backward path.
Then the entry in $\om'$ is moved to the left/top until $\om$ by
modified jeu de taquin. Let $P$ denote the backward path it takes.
Next comes the
application of (HC*2), where the
transformations on the entries of $H$ are performed, so that we
obtain a ``new" $H$. By induction, we
do not have to worry about the entries of the ``new" $H$ that are the
same as the corresponding entries in the ``old" $H$.
The only changes concerned the entries in column $j$ between rows
$i$ and $i'$.
Now we move to (HC*1) and determine candidate cells again. The entries
of the ``new" $H$ will obey the required lower bounds if and only if all these
``new" candidate cells are well-defined, i.e., if and only if they
lie {\it inside} the Ferrers diagram of $\la$. The only ``new"
candidate cells that we have to worry about are the ones between rows
$i$ and $i'$. However, it follows immediately from Lemma~HC* that
these ``new" candidate cells are located weakly to the left of the
backward path $P$, so,
at worst, these ``new" candidate cells are located {\it
on} $P$.
Consequently, they are still located inside the Ferrers diagram.
This establishes that $H$ is indeed a hook tabloid at any stage
during Algorithm~HC*.

\smallskip
What remains is to demonstrate that Algorithm~HC* and Algorithm~HC are
inverses of each other. We would like to show that the two
algorithms reverse each other step by step. That Algorithm~HC
reverses Algorithm~HC* is obvious. For the converse, we would need
the following: Consider the pair $(T,H)$ obtained after some loop
(HC1)-(HC2) during the performance of Algorithm~HC. In particular,
suppose that in the last application of (HC1) we moved from cell $\om$
to cell $\om'$ to obtain $T$. If we
now, with this pair $(T,H)$ and with distinguished cell equal to $\om$, 
jump into (HC*1) of Algorithm~HC*, then $\om'$ is the candidate cell
with first backward jeu de taquin path in the total order of
backwards paths. (It is a candidate cell, of course.) This fact is
demonstrated in Lemma~HC in Section~3.

\smallskip
Consequently, under the assumption of the truth of the lemmas in
Section~3, it is abundantly clear that Algorithms~HC and HC* are
well-defined and inverses of each other. This finishes the bijective
proof of (1.3).

\subhead 3. The crucial lemmas and their proofs\endsubhead

\proclaim{Lemma HC*}Let $\om$ be a cell of $\la$ which is not at the
bottom of some column, and let $\bom$ be its successor in the fixed
total order of the cells of $\la$. Furthermore,
at some stage during the execution of Algorithm~HC*,
let $P$ be the backward path of modified jeu de taquin which is performed in
(HC*1), connecting the candidate cell, $\om'$ say, with first
backward path in the total order of paths that was 
defined in Definition~P with cell
$\om$. Let $\bom'$ be a candidate cell in
the {\rm subsequent} execution of (HC*1). Then the backward path of modified
jeu de taquin, $\bP$ say, which starts in $\bom'$ reaches $\bom$.
Moreover, the backward paths $P$ and $\bP$
do not cross each other, and, in addition, horizontal pieces of $P$
and $\bP$ do not overlap.
\endproclaim
\demo{Proof}This claim follows, as we shall see, 
from the way the modifications of $H$ in (HC*2) are
defined and from a standard property of jeu de
taquin.

The statement of the lemma talks about two executions of (HC*1).
Suppose we are in the first of these two, i.e., we determine
candidate cells, of which $\om'$ is the one with first backward path in the
total order of backward paths.
The definition of the path ordering implies that any candidate cell
either is located in a lower row than $\om'$, and if so,
its backward path must enter the row which
contains $\om'$ weakly to the left of $\om'$, or, if not,
is located weakly to the left of the backward path $P$, but is allowed to be
on the backward path only on its horizontal pieces, with their rightmost cells
excluded (these are the cells on the backward path where it changes
direction from north to west), see Figure~6.a (there, possible candidate
cells are indicated by bold dots). Later, we will refer to these
candidate cells as the ``old" candidate cells.
\vskip10pt
\vbox{
$$
\Einheit.3cm
\PfadDicke{2pt}
\Pfad(0,-2),22\endPfad
\Pfad(0,-2),1111\endPfad
\Pfad(4,-5),222\endPfad
\Pfad(4,-5),111\endPfad
\Pfad(7,-7),22\endPfad
\Kreis(0,0)
\Kreis(0,-1)
\Kreis(4,-2)
\Kreis(4,-3)
\Kreis(4,-4)
\Kreis(7,-5)
\Kreis(7,-6)
\DickPunkt(0,-2)
\DickPunkt(1,-2)
\DickPunkt(2,-2)
\DickPunkt(3,-2)
\DickPunkt(0,-3)
\DickPunkt(1,-3)
\DickPunkt(2,-3)
\DickPunkt(3,-3)
\DickPunkt(0,-4)
\DickPunkt(1,-4)
\DickPunkt(2,-4)
\DickPunkt(3,-4)
\DickPunkt(0,-5)
\DickPunkt(1,-5)
\DickPunkt(2,-5)
\DickPunkt(3,-5)
\DickPunkt(4,-5)
\DickPunkt(5,-5)
\DickPunkt(6,-5)
\DickPunkt(0,-6)
\DickPunkt(1,-6)
\DickPunkt(2,-6)
\DickPunkt(3,-6)
\DickPunkt(4,-6)
\DickPunkt(5,-6)
\DickPunkt(6,-6)
\DickPunkt(0,-7)
\DickPunkt(1,-7)
\DickPunkt(2,-7)
\DickPunkt(3,-7)
\DickPunkt(4,-7)
\DickPunkt(5,-7)
\DickPunkt(6,-7)
\DickPunkt(7,-7)
\DickPunkt(0,-8)
\DickPunkt(1,-8)
\DickPunkt(2,-8)
\DickPunkt(3,-8)
\DickPunkt(4,-8)
\DickPunkt(5,-8)
\DickPunkt(6,-8)
\DickPunkt(7,-8)
\DickPunkt(8,-8)
\DickPunkt(9,-8)
\DickPunkt(10,-8)
\DickPunkt(0,-9)
\DickPunkt(1,-9)
\DickPunkt(2,-9)
\DickPunkt(3,-9)
\DickPunkt(4,-9)
\DickPunkt(5,-9)
\DickPunkt(6,-9)
\DickPunkt(7,-9)
\DickPunkt(8,-9)
\DickPunkt(9,-9)
\DickPunkt(10,-9)
\Label\o{\om}(0,0)
\Label\ro{\kern7pt\om'}(7,-7)
\Label\r{\hskip2pt\hbox to 1cm{\leaders\hbox to .3cm{\hss.\hss}\hfill}}(11,-8)
\Label\r{\hskip2pt\hbox to 1cm{\leaders\hbox to .3cm{\hss.\hss}\hfill}}(11,-9)
\Label\u{\hbox to 4.5cm{\leaders\hbox to .3cm{\hss.\hss}\hfill}}(7,-9)
\hbox{\hskip5.5cm}
\Pfad(0,-2),22\endPfad
\Pfad(0,-2),1111\endPfad
\Pfad(4,-5),222\endPfad
\Pfad(4,-5),111\endPfad
\Pfad(7,-7),22\endPfad
\Kreis(0,0)
\DickPunkt(0,-1)
\Kreis(4,-2)
\DickPunkt(4,-3)
\DickPunkt(4,-4)
\Kreis(7,-5)
\DickPunkt(7,-6)
\DickPunkt(0,-2)
\Kreis(1,-2)
\Kreis(2,-2)
\Kreis(3,-2)
\DickPunkt(0,-3)
\DickPunkt(1,-3)
\DickPunkt(2,-3)
\DickPunkt(3,-3)
\DickPunkt(0,-4)
\DickPunkt(1,-4)
\DickPunkt(2,-4)
\DickPunkt(3,-4)
\DickPunkt(0,-5)
\DickPunkt(1,-5)
\DickPunkt(2,-5)
\DickPunkt(3,-5)
\DickPunkt(4,-5)
\Kreis(5,-5)
\Kreis(6,-5)
\DickPunkt(0,-6)
\DickPunkt(1,-6)
\DickPunkt(2,-6)
\DickPunkt(3,-6)
\DickPunkt(4,-6)
\DickPunkt(5,-6)
\DickPunkt(6,-6)
\DickPunkt(0,-7)
\DickPunkt(1,-7)
\DickPunkt(2,-7)
\DickPunkt(3,-7)
\DickPunkt(4,-7)
\DickPunkt(5,-7)
\DickPunkt(6,-7)
\DickPunkt(7,-7)
\DickPunkt(0,-8)
\DickPunkt(1,-8)
\DickPunkt(2,-8)
\DickPunkt(3,-8)
\DickPunkt(4,-8)
\DickPunkt(5,-8)
\DickPunkt(6,-8)
\DickPunkt(7,-8)
\DickPunkt(8,-8)
\DickPunkt(9,-8)
\DickPunkt(10,-8)
\DickPunkt(0,-9)
\DickPunkt(1,-9)
\DickPunkt(2,-9)
\DickPunkt(3,-9)
\DickPunkt(4,-9)
\DickPunkt(5,-9)
\DickPunkt(6,-9)
\DickPunkt(7,-9)
\DickPunkt(8,-9)
\DickPunkt(9,-9)
\DickPunkt(10,-9)
\Label\o{\om}(0,0)
\Label\ro{\kern7pt\om'}(7,-7)
\Label\r{\hskip2pt\hbox to 1cm{\leaders\hbox to .3cm{\hss.\hss}\hfill}}(11,-8)
\Label\r{\hskip2pt\hbox to 1cm{\leaders\hbox to .3cm{\hss.\hss}\hfill}}(11,-9)
\Label\u{\hbox to 4.5cm{\leaders\hbox to .3cm{\hss.\hss}\hfill}}(7,-9)
\hskip3.9cm
$$
\centerline{\eightpoint a. possible ``old" candidate cells\hskip1.5cm
b. possible ``new" candidate cells}
\vskip10pt
\centerline{\eightpoint Figure 6}
}
\vskip10pt

Next comes the
application of (HC*2), where the
transformations on the entries of $H$ are performed. 
Then we move to (HC*1) and determine candidate cells again. Because
of the definition of the transformations on $H$ in (HC*2), these ``new"
candidate cells are related to the ``old" ones as follows: Those which
were located in a lower row than $\om'$ are still
candidate cells. All the other ``old" candidate cells, with the
exception of $\om'$, moved one unit to the right and
one unit down to form the corresponding ``new" candidate
cells. Therefore, a ``new" candidate cell is either located 
in a lower row than $\om'$, and if so,
its backward path must enter the row which
contains $\om'$ weakly to the left of $\om'$, or, if not,
is located weakly to the left of the backward path $P$, but is allowed to be
on the backward path only on its vertical pieces, with their topmost cells
excluded (these are again the cells on the backward path where it changes
direction from north to west), see Figure~6.b (possible ``new" candidate
cells are indicated by bold dots).

The assertions of the lemma now follow from the following property of
our modified jeu de taquin: If the second backward 
path is below of the first backward path somewhere,
then it has to stay below from thereon. To make this precise,
suppose that during the first execution of (HC*1) the special entry, 
$s_1$ say, went to the left
by the move (2.6)-(2.7), see the left half of Figure~7. (The arrows mark the
direction of move of the special entry.)
\vskip10pt
\vbox{
$$
\smatrix \format \sa\c\s\c\se\\
\hlinefor5\\
&y&&s_1&\\
\hlinefor5\\
&z&&*&\\
\hlinefor5
\endsmatrix
\overset \text{during}\to{\overset \text{first (HC*1)}\to\longrightarrow}
\smatrix \format \sa\c\s\c\se\\
\hlinefor5\\
&\hbox to5pt{$*\kern-4pt\leftarrow$\hss}&&y-1&\\
\hlinefor5\\
&\,z\,&&*&\\
\hlinefor5
\endsmatrix
\hskip2cm
\smatrix \format \sa\c\s\c\se\\
\hlinefor5\\
&*&&y-1&\\
\hlinefor5\\
&z&&s_2&\\
\hlinefor5
\endsmatrix
\overset \text{during}\to {\overset \text{second (HC*1)}\to\longrightarrow}
\smatrix \format \sa\c\s\c\se\\
\hlinefor5\\
&\,*\,&&y-1&\\
\hlinefor5\\
&\hbox to5pt{$*\kern-4pt\leftarrow$\hss}&&z-1&\\
\hlinefor5
\endsmatrix
$$
\centerline{\eightpoint Figure 7}
}
\vskip10pt
\noindent
Since columns are strictly increasing, we have $y<z$. Suppose that
during the second execution of (HC*1) we reach the cell neighbouring 
$z$ and $y-1$ with
a special entry $s_2$, see the right half of Figure~7. Then the
definition of (JT*) forces us to move left in the next step,
i.e., to apply the move (2.6)-(2.7). 

This concludes the proof of the lemma.\quad \quad \qed
\enddemo

\proclaim{Lemma OC*}The output of Algorithm~HC*
is a content tabloid.
\endproclaim
\demo{Proof} Let the output of Algorithm~HC* be the array $C$.
We have to prove that the entry $C_\rh$ of $C$ in cell
$\rh$ satisfies 
$$1-c_\rh\le C_\rh\le b.\tag3.1$$
Observe first that once an entry is moved into the distinguished cell 
during the
execution of (HC*1)-(HC*2), it stays there for the rest 
of Algorithm~HC*, and will therefore be the corresponding
entry in the output $C$. 

It is easy to check that the moves (2.6)-(2.7) and (2.6)-(2.8)
guarantee that after each loop (HC*1)-(HC*2) 
the portion of the array $T$, obtained
after that loop, consisting of all cells which are after
the distinguished cell in the fixed total order of the cells, is a {\it skew
tableau\/}, i.e., the entries are weakly increasing along rows and
the entries are strictly increasing along columns in that portion.
Besides, if the distinguished cell is the bottommost cell in column $j_0$,
say, then any
entry in the first row and to the right of the distinguished cell is at least $1-j_0$.
This follows, by induction on $j_0$, from the definition of $T$ in
(HC*0) as a tableau with entries at least 1 and the fact that it 
is only the move (2.6)-(2.7) where something (namely 1) is subtracted
from an entry. If we combine both facts, we obtain that
if the distinguished cell is the bottommost cell in column $j_0$, then any
entry in row $i_0$ and to the right of the distinguished cell 
is at least $i_0-j_0$.

Now, suppose that we are in (HC*1) and are about to move an entry, $s$
say, into the distinguished cell,
$\om=(i,j)$ say. If $i=1$, i.e., if $\om$ is in the first row, then the
observation of the preceding paragraph (with $j_0=j-1$) 
immediately implies that this
entry $s$ is at least $1-(j-1)=1-c_{(1,j)}$, which establishes the
lemma in that case. If $\om$ is not in the first row, then let
$\widehat\om=(i-1,j)$ be the preceding cell in the fixed total order of the
cells, let $\widehat P$ be the backward path into $\widehat\om$ of the preceding
application of (HC*1), and let $\widehat\om'$ be the cell where 
$\widehat P$ starts. Lemma~HC* implies that the entry $s$ comes from a
cell, $(\widehat i,\widehat j)$ say, which is located in a lower row than $\widehat\om'$, or, if not,
is located weakly to the left of the backward path $\widehat P$, but is allowed to be
on the backward path only on its vertical pieces, with their topmost cells
excluded. Furthermore, that cell $(\widehat i,\widehat j)$ is not on a
horizontal piece of any other ``previous" backward path which terminated in column $j$. 
Consequently, again by the observation of the preceding paragraph
(with $i_0=\widehat i$ and $j_0=j-1$) we have $s\ge \widehat i-(j-1)\ge
i-(j-1)=1-c_{(i,j)}$, which is exactly the assertion of the lemma.\quad 
\quad \qed
\enddemo

\proclaim{Lemma HC}Suppose that we are at the end of a loop
(HC1)-(HC2) during the execution of Algorithm~HC. In particular,
suppose that in the last application of (HC1) we moved from cell $\om$
to cell $\om'$, and that we obtained the pair $(T,H)$ in (HC2). If we
now, with this pair $(T,H)$ and with distinguished cell equal to $\om$, 
jump into (HC*1) of Algorithm~HC*, then $\om'$ is the candidate cell
with first backward jeu de taquin path in the total order of
backwards paths that was
defined in Definition~P in Section~2.
\endproclaim
\demo{Proof} We prove the assertion of the lemma by {\it reverse} induction on the
number of the cell $\om$ in the fixed total order of the cells.

There is nothing to prove if $\om$ is the bottommost
cell in some column. In particular, the assertion is true if $\om$ is
the very last cell in the total order of cells, i.e., the bottommost
cell in the last column of $\la$. This allows to start the induction.

Now, let us assume that $\om$ is not the
bottommost cell in some column. Then the successor of $\om$, $\bom$
say, in the fixed total order of cells is in the same column
immediately below of
$\om$. Let $(\bT,\bH)$ be the outcome of the preceding loop
(HC1)-(HC2), i.e., we obtain $(T,H)$ by applying (HC1)-(HC2) to
$(\bT,\bH)$ with distinguished cell equal to $\om$. On the other
hand, $\bT$ was obtained as the outcome of the previous execution of
(HC1), when some entry was moved by (JT) from $\bom$ to $\bom'$, say.
Let $\bP$ be this foward path.
By induction, we may assume that, when we jump into (HC*1) with
$(\bT,\bH)$ and distinguished cell $\bom$, the backward path of $\bom'$ is the
first among all backward paths of candidate cells. From now on, we will
refer to these candidate cells as the ``old" candidate cells.

In particular, as we already observed in the proof of Lemma~HC*, this
says that any ``old" candidate cell
either is located in a lower row than $\bom'$, and if so,
its backward path must enter the row which
contains $\bom'$ weakly to the left of $\bom'$, or, if not,
is located weakly to the left of the forward path $\bP$, but is allowed to be
on the forward path only on its horizontal pieces, with their rightmost cells
excluded, see Figure~8.a. (There, possible ``old" candidate cells are
indicated by bold dots.) 
\vskip10pt
\vbox{
$$
\Einheit.3cm
\PfadDicke{2pt}
\Pfad(0,-2),2\endPfad
\Pfad(0,-2),1111\endPfad
\Pfad(4,-4),22\endPfad
\Pfad(4,-4),11\endPfad
\Kreis(0,0)
\Kreis(0,-1)
\Kreis(4,-2)
\Kreis(4,-3)
\Kreis(4,-4)
\DickPunkt(0,-2)
\DickPunkt(1,-2)
\DickPunkt(2,-2)
\DickPunkt(3,-2)
\DickPunkt(0,-3)
\DickPunkt(1,-3)
\DickPunkt(2,-3)
\DickPunkt(3,-3)
\DickPunkt(0,-4)
\DickPunkt(1,-4)
\DickPunkt(2,-4)
\DickPunkt(3,-4)
\DickPunkt(4,-4)
\DickPunkt(5,-4)
\DickPunkt(6,-4)
\DickPunkt(0,-5)
\DickPunkt(1,-5)
\DickPunkt(2,-5)
\DickPunkt(3,-5)
\DickPunkt(4,-5)
\DickPunkt(5,-5)
\DickPunkt(6,-5)
\DickPunkt(7,-5)
\DickPunkt(8,-5)
\DickPunkt(9,-5)
\DickPunkt(10,-5)
\DickPunkt(11,-5)
\DickPunkt(12,-5)
\DickPunkt(13,-5)
\DickPunkt(14,-5)
\DickPunkt(0,-6)
\DickPunkt(1,-6)
\DickPunkt(2,-6)
\DickPunkt(3,-6)
\DickPunkt(4,-6)
\DickPunkt(5,-6)
\DickPunkt(6,-6)
\DickPunkt(7,-6)
\DickPunkt(8,-6)
\DickPunkt(9,-6)
\DickPunkt(10,-6)
\DickPunkt(10,-6)
\DickPunkt(11,-6)
\DickPunkt(12,-6)
\DickPunkt(13,-6)
\DickPunkt(14,-6)
\DickPunkt(0,-7)
\DickPunkt(1,-7)
\DickPunkt(2,-7)
\DickPunkt(3,-7)
\DickPunkt(4,-7)
\DickPunkt(5,-7)
\DickPunkt(6,-7)
\DickPunkt(7,-7)
\DickPunkt(8,-7)
\DickPunkt(9,-7)
\DickPunkt(10,-7)
\DickPunkt(11,-7)
\DickPunkt(12,-7)
\DickPunkt(13,-7)
\DickPunkt(14,-7)
\DickPunkt(0,-8)
\DickPunkt(1,-8)
\DickPunkt(2,-8)
\DickPunkt(3,-8)
\DickPunkt(4,-8)
\DickPunkt(5,-8)
\DickPunkt(6,-8)
\DickPunkt(7,-8)
\DickPunkt(8,-8)
\DickPunkt(9,-8)
\DickPunkt(10,-8)
\DickPunkt(11,-8)
\DickPunkt(12,-8)
\DickPunkt(13,-8)
\DickPunkt(14,-8)
\DickPunkt(0,-9)
\DickPunkt(1,-9)
\DickPunkt(2,-9)
\DickPunkt(3,-9)
\DickPunkt(4,-9)
\DickPunkt(5,-9)
\DickPunkt(6,-9)
\DickPunkt(7,-9)
\DickPunkt(8,-9)
\DickPunkt(9,-9)
\DickPunkt(10,-9)
\DickPunkt(11,-9)
\DickPunkt(12,-9)
\DickPunkt(13,-9)
\DickPunkt(14,-9)
\DickPunkt(0,-10)
\DickPunkt(1,-10)
\DickPunkt(2,-10)
\DickPunkt(3,-10)
\DickPunkt(4,-10)
\DickPunkt(5,-10)
\DickPunkt(6,-10)
\DickPunkt(7,-10)
\DickPunkt(8,-10)
\DickPunkt(9,-10)
\DickPunkt(10,-10)
\DickPunkt(11,-10)
\DickPunkt(12,-10)
\DickPunkt(13,-10)
\DickPunkt(14,-10)
\Label\o{\bP}(2,-2)
\Label\o{\om}(0,0)
\Label\l{\bom}(0,-1)
\Label\r{\bom'}(6,-4)
\Label\r{\hskip2pt\hbox to 1cm{\leaders\hbox to
.3cm{\hss.\hss}\hfill}}(15,-5)
\Label\r{\hskip2pt\hbox to 1cm{\leaders\hbox to
.3cm{\hss.\hss}\hfill}}(15,-6)
\Label\r{\hskip2pt\hbox to 1cm{\leaders\hbox to
.3cm{\hss.\hss}\hfill}}(15,-7)
\Label\r{\hskip2pt\hbox to 1cm{\leaders\hbox to
.3cm{\hss.\hss}\hfill}}(15,-8)
\Label\r{\hskip2pt\hbox to 1cm{\leaders\hbox to
.3cm{\hss.\hss}\hfill}}(15,-9)
\Label\r{\hskip2pt\hbox to 1cm{\leaders\hbox to
.3cm{\hss.\hss}\hfill}}(15,-10)
\Label\u{\hbox to 5.7cm{\leaders\hbox to
.3cm{\hss.\hss}\hfill}}(9,-10)
\hbox{\hskip6.5cm}
\Pfad(0,-2),2\endPfad
\Pfad(0,-2),1111\endPfad
\Pfad(4,-4),22\endPfad
\Pfad(4,-4),11\endPfad
\Kreis(0,-2)
\Kreis(1,-2)
\Kreis(2,-2)
\Kreis(3,-2)
\Kreis(4,-4)
\Kreis(5,-4)
\DickPunkt(0,0)
\DickPunkt(1,0)
\DickPunkt(2,0)
\DickPunkt(3,0)
\DickPunkt(4,0)
\DickPunkt(5,0)
\DickPunkt(6,0)
\DickPunkt(7,0)
\DickPunkt(8,0)
\DickPunkt(9,0)
\DickPunkt(10,0)
\DickPunkt(11,0)
\DickPunkt(12,0)
\DickPunkt(13,0)
\DickPunkt(14,0)
\DickPunkt(0,-1)
\DickPunkt(1,-1)
\DickPunkt(2,-1)
\DickPunkt(3,-1)
\DickPunkt(4,-1)
\DickPunkt(5,-1)
\DickPunkt(6,-1)
\DickPunkt(7,-1)
\DickPunkt(8,-1)
\DickPunkt(9,-1)
\DickPunkt(10,-1)
\DickPunkt(11,-1)
\DickPunkt(12,-1)
\DickPunkt(13,-1)
\DickPunkt(14,-1)
\DickPunkt(4,-2)
\DickPunkt(5,-2)
\DickPunkt(6,-2)
\DickPunkt(7,-2)
\DickPunkt(8,-2)
\DickPunkt(9,-2)
\DickPunkt(10,-2)
\DickPunkt(11,-2)
\DickPunkt(12,-2)
\DickPunkt(13,-2)
\DickPunkt(14,-2)
\DickPunkt(4,-3)
\DickPunkt(5,-3)
\DickPunkt(6,-3)
\DickPunkt(7,-3)
\DickPunkt(8,-3)
\DickPunkt(9,-3)
\DickPunkt(10,-3)
\DickPunkt(11,-3)
\DickPunkt(12,-3)
\DickPunkt(13,-3)
\DickPunkt(14,-3)
\DickPunkt(6,-4)
\DickPunkt(7,-4)
\DickPunkt(8,-4)
\DickPunkt(9,-4)
\DickPunkt(10,-4)
\DickPunkt(11,-4)
\DickPunkt(12,-4)
\DickPunkt(13,-4)
\DickPunkt(14,-4)
\DickPunkt(6,-5)
\DickPunkt(7,-5)
\DickPunkt(8,-5)
\DickPunkt(9,-5)
\DickPunkt(10,-5)
\DickPunkt(11,-5)
\DickPunkt(12,-5)
\DickPunkt(13,-5)
\DickPunkt(14,-5)
\DickPunkt(6,-6)
\DickPunkt(7,-6)
\DickPunkt(8,-6)
\DickPunkt(9,-6)
\DickPunkt(10,-6)
\DickPunkt(11,-6)
\DickPunkt(12,-6)
\DickPunkt(13,-6)
\DickPunkt(14,-6)
\DickPunkt(6,-7)
\DickPunkt(7,-7)
\DickPunkt(8,-7)
\DickPunkt(9,-7)
\DickPunkt(10,-7)
\DickPunkt(11,-7)
\DickPunkt(12,-7)
\DickPunkt(13,-7)
\DickPunkt(14,-7)
\DickPunkt(6,-8)
\DickPunkt(7,-8)
\DickPunkt(8,-8)
\DickPunkt(9,-8)
\DickPunkt(10,-8)
\DickPunkt(11,-8)
\DickPunkt(12,-8)
\DickPunkt(13,-8)
\DickPunkt(14,-8)
\DickPunkt(6,-9)
\DickPunkt(7,-9)
\DickPunkt(8,-9)
\DickPunkt(9,-9)
\DickPunkt(10,-9)
\DickPunkt(11,-9)
\DickPunkt(12,-9)
\DickPunkt(13,-9)
\DickPunkt(14,-9)
\DickPunkt(6,-10)
\DickPunkt(7,-10)
\DickPunkt(8,-10)
\DickPunkt(9,-10)
\DickPunkt(10,-10)
\DickPunkt(11,-10)
\DickPunkt(12,-10)
\DickPunkt(13,-10)
\DickPunkt(14,-10)
\Label\u{\bP}(2,-2)
\Label\o{\om}(0,0)
\Label\l{\bom}(0,-1)
\Label\l{\bom'}(6,-5)
\Label\r{\hskip2pt\hbox to 1cm{\leaders\hbox to
.3cm{\hss.\hss}\hfill}}(15,0)
\Label\r{\hskip2pt\hbox to 1cm{\leaders\hbox to
.3cm{\hss.\hss}\hfill}}(15,-1)
\Label\r{\hskip2pt\hbox to 1cm{\leaders\hbox to
.3cm{\hss.\hss}\hfill}}(15,-2)
\Label\r{\hskip2pt\hbox to 1cm{\leaders\hbox to
.3cm{\hss.\hss}\hfill}}(15,-3)
\Label\r{\hskip2pt\hbox to 1cm{\leaders\hbox to
.3cm{\hss.\hss}\hfill}}(15,-4)
\Label\r{\hskip2pt\hbox to 1cm{\leaders\hbox to
.3cm{\hss.\hss}\hfill}}(15,-5)
\Label\r{\hskip2pt\hbox to 1cm{\leaders\hbox to
.3cm{\hss.\hss}\hfill}}(15,-6)
\Label\r{\hskip2pt\hbox to 1cm{\leaders\hbox to
.3cm{\hss.\hss}\hfill}}(15,-7)
\Label\r{\hskip2pt\hbox to 1cm{\leaders\hbox to
.3cm{\hss.\hss}\hfill}}(15,-8)
\Label\r{\hskip2pt\hbox to 1cm{\leaders\hbox to
.3cm{\hss.\hss}\hfill}}(15,-9)
\Label\r{\hskip2pt\hbox to 1cm{\leaders\hbox to
.3cm{\hss.\hss}\hfill}}(15,-10)
\Label\u{\hbox to 3.9cm{\leaders\hbox to
.3cm{\hss.\hss}\hfill}}(12,-10)
\hbox{\hskip5.1cm}
$$
\centerline{\eightpoint a. possible ``old" candidate cells\hskip1.5cm
b. possible territory for $P$}
\vskip10pt
\centerline{\eightpoint Figure 8}
}
\vskip10pt

Next, in the execution of Algorithm~HC, comes the
application of (HC1)-(HC2) to obtain $(T,H)$ out of $(\bT,\bH)$.
Let $P$ be the forward path from $\om$ to $\om'$ performed in (HC1). 
We claim that $P$ has to be ``above" $\bP$ always. To make this
precise, suppose that during the first execution of (HC1) the special entry, 
$s_1$ say, on its way along $\bP$ went to the right
by the move (2.2)-(2.3), see the left half of Figure~9. (The arrows mark the
direction of move of the special entry.)
\vskip10pt
\vbox{
$$
\smatrix \format \sa\c\s\c\se\\
\hlinefor5\\
&*&&y&\\
\hlinefor5\\
&s_1&&z&\\
\hlinefor5
\endsmatrix
\overset \text{during}\to{\overset \text{first (HC1)}\to\longrightarrow}
\smatrix \format \sa\c\s\c\se\\
\hlinefor5\\
&*&&\,y\,&\\
\hlinefor5\\
&z+1&&\hbox to5pt{\hss$\rightarrow\kern-4pt*$}&\\
\hlinefor5
\endsmatrix
\hskip2cm
\smatrix \format \sa\c\s\c\se\\
\hlinefor5\\
&s_2&&y&\\
\hlinefor5\\
&z+1&&*&\\
\hlinefor5
\endsmatrix
\overset \text{during}\to {\overset \text{second (HC1)}\to\longrightarrow}
\smatrix \format \sa\c\s\c\se\\
\hlinefor5\\
&y+1&&\hbox to5pt{\hss$\rightarrow\kern-4pt*$}&\\
\hlinefor5\\
&z+1&&\,*\,&\\
\hlinefor5
\endsmatrix
$$
\centerline{\eightpoint Figure 9}
}
\vskip10pt
\noindent
Since columns are strictly increasing, we have $y<z$. Suppose that
during the second execution of (HC1), when moving along $P$, 
we reach the cell neighbouring 
$z+1$ and $y$ with
a special entry $s_2$, see the right half of Figure~9. Then the
definition of (JT) forces us to move right in the next step,
i.e., to apply the move (2.2)-(2.3). As a corollary, we obtain that
$P$ must be located in the region which is above and to the right of
$\bP$, with no overlap of horizontal pieces of $P$ and $\bP$, 
as indicated in Figure~8.b. (There, the cells through which $P$ may
run are indicated by bold dots.)

Furthermore, we claim that $P$ has to be ``above" the backward path,
of {\it any} ``old" candidate cell. 
To make this precise, suppose that we jump into (HC*1) with $(\bT,\bH)$
and distinguished cell $\bom$ and determine the backward path, $Q$
say, of an ``old" candidate cell, $\ze$ say. Let us further suppose
that, while determining $Q$, the special entry,
$s_1$ say, would go to the left
by the move (2.6)-(2.7), see the left half of Figure~10. 
(The arrows mark the
direction of (potential) move of the special entry.)
\vskip10pt
\vbox{
$$
\overset \text{during}\to{\overset \text{(HC*1)}\to{}}
\smatrix \format \sa\c\s\c\se\\
\hlinefor5\\
&*&&\,y\,&\\
\hlinefor5\\
&\hbox to5pt{$z\kern-4pt\leftarrow$\hss}&&s_1&\\
\hlinefor5
\endsmatrix
\hskip2cm
\smatrix \format \sa\c\s\c\se\\
\hlinefor5\\
&s_2&&y&\\
\hlinefor5\\
&z&&*&\\
\hlinefor5
\endsmatrix
\overset \text{during}\to {\overset \text{(HC1)}\to\longrightarrow}
\smatrix \format \sa\c\s\c\se\\
\hlinefor5\\
&y+1&&\hbox to5pt{\hss$\rightarrow\kern-4pt*$}&\\
\hlinefor5\\
&z&&\,s_1\,&\\
\hlinefor5
\endsmatrix
$$
\centerline{\eightpoint Figure 10}
}
\vskip10pt
\noindent
If this is the case then we must have $z-1>y$. Now suppose that we
apply (HC1)-(HC2) to $(\bT,\bH)$, with distinguished cell $\bom$, 
thus obtaining $(T,H)$, and now
apply (HC1), with distinguished cell $\om$. Suppose further that
during the execution of (HC1), when moving along $P$, 
we reach the cell neighbouring 
$z$ and $y$ with
a special entry $s_2$, see the right half of Figure~10. Then the
definition of (JT) forces us to move to the right in the next step,
i.e., to apply the move (2.2)-(2.3). 

Combining the above with the already observed fact that backward
paths of ``old" candidate cells in a lower row than $\bom'$ must 
enter the row of $\bom'$ weakly to the left of $\bom'$, we obtain that
``old" candidate cells are located in one of the following three
regions (see Figure~11.a; there, possible ``old" candidate cells are
indicated by bold dots):
\roster
\item "(I)" In a row strictly below $\om'$; then the backward path of
the ``old" candidate cell must enter the row of $\om'$ weakly to the
left of $\om'$.
\item "(II)" In a row strictly below $\bom'$ and weakly above $\om'$,
weakly to the left of $P$, but on $P$ only on its vertical
pieces, with their topmost cells excluded.
\item "(III)" In a row weakly above $\bom'$ and weakly to the left of
$\bP$, but on $\bP$ only on its horizontal pieces, with their rightmost cells
excluded.
\endroster
\vskip10pt
\vbox{
$$
\Einheit.3cm
\PfadDicke{2pt}
\Pfad(0,-2),2\endPfad
\Pfad(0,-2),1111\endPfad
\Pfad(4,-4),22\endPfad
\Pfad(4,-4),11\endPfad
\hbox{\hskip1pt}
\Pfad(0,0),1111\endPfad
\Pfad(4,-3),222\endPfad
\Pfad(4,-3),1111\endPfad
\Pfad(8,-5),22\endPfad
\Pfad(8,-5),111\endPfad
\Pfad(11,-7),22\endPfad
\Pfad(11,-7),11\endPfad
\Pfad(13,-9),22\endPfad
\hbox{\hskip-1pt}
\Kreis(0,0)
\Kreis(0,-1)
\Kreis(4,-2)
\Kreis(4,-3)
\Kreis(4,-4)
\Kreis(9,-5)
\Kreis(10,-5)
\Kreis(11,-5)
\Kreis(12,-7)
\Kreis(13,-7)
\DickPunkt(0,-2)
\DickPunkt(1,-2)
\DickPunkt(2,-2)
\DickPunkt(3,-2)
\DickPunkt(0,-3)
\DickPunkt(1,-3)
\DickPunkt(2,-3)
\DickPunkt(3,-3)
\DickPunkt(0,-4)
\DickPunkt(1,-4)
\DickPunkt(2,-4)
\DickPunkt(3,-4)
\DickPunkt(4,-4)
\DickPunkt(5,-4)
\DickPunkt(6,-4)
\DickPunkt(0,-5)
\DickPunkt(1,-5)
\DickPunkt(2,-5)
\DickPunkt(3,-5)
\DickPunkt(4,-5)
\DickPunkt(5,-5)
\DickPunkt(6,-5)
\DickPunkt(7,-5)
\DickPunkt(8,-5)
\DickPunkt(0,-6)
\DickPunkt(1,-6)
\DickPunkt(2,-6)
\DickPunkt(3,-6)
\DickPunkt(4,-6)
\DickPunkt(5,-6)
\DickPunkt(6,-6)
\DickPunkt(7,-6)
\DickPunkt(8,-6)
\DickPunkt(9,-6)
\DickPunkt(10,-6)
\DickPunkt(10,-6)
\DickPunkt(11,-6)
\DickPunkt(0,-7)
\DickPunkt(1,-7)
\DickPunkt(2,-7)
\DickPunkt(3,-7)
\DickPunkt(4,-7)
\DickPunkt(5,-7)
\DickPunkt(6,-7)
\DickPunkt(7,-7)
\DickPunkt(8,-7)
\DickPunkt(9,-7)
\DickPunkt(10,-7)
\DickPunkt(11,-7)
\DickPunkt(0,-8)
\DickPunkt(1,-8)
\DickPunkt(2,-8)
\DickPunkt(3,-8)
\DickPunkt(4,-8)
\DickPunkt(5,-8)
\DickPunkt(6,-8)
\DickPunkt(7,-8)
\DickPunkt(8,-8)
\DickPunkt(9,-8)
\DickPunkt(10,-8)
\DickPunkt(11,-8)
\DickPunkt(12,-8)
\DickPunkt(13,-8)
\DickPunkt(0,-9)
\DickPunkt(1,-9)
\DickPunkt(2,-9)
\DickPunkt(3,-9)
\DickPunkt(4,-9)
\DickPunkt(5,-9)
\DickPunkt(6,-9)
\DickPunkt(7,-9)
\DickPunkt(8,-9)
\DickPunkt(9,-9)
\DickPunkt(10,-9)
\DickPunkt(11,-9)
\DickPunkt(12,-9)
\DickPunkt(13,-9)
\DickPunkt(0,-10)
\DickPunkt(1,-10)
\DickPunkt(2,-10)
\DickPunkt(3,-10)
\DickPunkt(4,-10)
\DickPunkt(5,-10)
\DickPunkt(6,-10)
\DickPunkt(7,-10)
\DickPunkt(8,-10)
\DickPunkt(9,-10)
\DickPunkt(10,-10)
\DickPunkt(11,-10)
\DickPunkt(12,-10)
\DickPunkt(13,-10)
\DickPunkt(14,-10)
\Label\o{P}(8,-3)
\Label\o{\bP}(2,-2)
\Label\o{\om}(0,0)
\Label\l{\bom}(0,-1)
\Label\r{\bom'}(6,-4)
\Label\ro{\kern7pt\om'}(13,-9)
\Label\r{\hskip2pt\hbox to 1cm{\leaders\hbox to
.3cm{\hss.\hss}\hfill}}(15,-10)
\Label\u{\hbox to 5.7cm{\leaders\hbox to
.3cm{\hss.\hss}\hfill}}(9,-10)
\hbox{\hskip6.5cm}
\Pfad(0,-2),2\endPfad
\Pfad(0,-2),1111\endPfad
\Pfad(4,-4),22\endPfad
\Pfad(4,-4),11\endPfad
\hbox{\hskip1pt}
\Pfad(0,0),1111\endPfad
\Pfad(4,-3),222\endPfad
\Pfad(4,-3),1111\endPfad
\Pfad(8,-5),22\endPfad
\Pfad(8,-5),111\endPfad
\Pfad(11,-7),22\endPfad
\Pfad(11,-7),11\endPfad
\Pfad(13,-9),22\endPfad
\hbox{\hskip-1pt}
\Kreis(4,0)
\Kreis(4,-1)
\Kreis(4,-2)
\Kreis(8,-3)
\Kreis(8,-4)
\Kreis(11,-5)
\Kreis(11,-6)
\Kreis(13,-7)
\Kreis(13,-8)
\DickPunkt(0,0)
\DickPunkt(1,0)
\DickPunkt(2,0)
\DickPunkt(3,0)
\DickPunkt(0,-1)
\DickPunkt(1,-1)
\DickPunkt(2,-1)
\DickPunkt(3,-1)
\DickPunkt(0,-2)
\DickPunkt(1,-2)
\DickPunkt(2,-2)
\DickPunkt(3,-2)
\DickPunkt(0,-3)
\DickPunkt(1,-3)
\DickPunkt(2,-3)
\DickPunkt(3,-3)
\DickPunkt(4,-3)
\DickPunkt(5,-3)
\DickPunkt(6,-3)
\DickPunkt(7,-3)
\DickPunkt(0,-4)
\DickPunkt(1,-4)
\DickPunkt(2,-4)
\DickPunkt(3,-4)
\DickPunkt(4,-4)
\DickPunkt(5,-4)
\DickPunkt(6,-4)
\DickPunkt(7,-4)
\DickPunkt(0,-5)
\DickPunkt(1,-5)
\DickPunkt(2,-5)
\DickPunkt(3,-5)
\DickPunkt(4,-5)
\DickPunkt(5,-5)
\DickPunkt(6,-5)
\DickPunkt(7,-5)
\DickPunkt(8,-5)
\DickPunkt(9,-5)
\DickPunkt(10,-5)
\DickPunkt(0,-6)
\DickPunkt(1,-6)
\DickPunkt(2,-6)
\DickPunkt(3,-6)
\DickPunkt(4,-6)
\DickPunkt(5,-6)
\DickPunkt(6,-6)
\DickPunkt(7,-6)
\DickPunkt(8,-6)
\DickPunkt(9,-6)
\DickPunkt(10,-6)
\DickPunkt(10,-6)
\DickPunkt(0,-7)
\DickPunkt(1,-7)
\DickPunkt(2,-7)
\DickPunkt(3,-7)
\DickPunkt(4,-7)
\DickPunkt(5,-7)
\DickPunkt(6,-7)
\DickPunkt(7,-7)
\DickPunkt(8,-7)
\DickPunkt(9,-7)
\DickPunkt(10,-7)
\DickPunkt(11,-7)
\DickPunkt(12,-7)
\DickPunkt(0,-8)
\DickPunkt(1,-8)
\DickPunkt(2,-8)
\DickPunkt(3,-8)
\DickPunkt(4,-8)
\DickPunkt(5,-8)
\DickPunkt(6,-8)
\DickPunkt(7,-8)
\DickPunkt(8,-8)
\DickPunkt(9,-8)
\DickPunkt(10,-8)
\DickPunkt(11,-8)
\DickPunkt(12,-8)
\DickPunkt(0,-9)
\DickPunkt(1,-9)
\DickPunkt(2,-9)
\DickPunkt(3,-9)
\DickPunkt(4,-9)
\DickPunkt(5,-9)
\DickPunkt(6,-9)
\DickPunkt(7,-9)
\DickPunkt(8,-9)
\DickPunkt(9,-9)
\DickPunkt(10,-9)
\DickPunkt(11,-9)
\DickPunkt(12,-9)
\DickPunkt(13,-9)
\DickPunkt(0,-10)
\DickPunkt(1,-10)
\DickPunkt(2,-10)
\DickPunkt(3,-10)
\DickPunkt(4,-10)
\DickPunkt(5,-10)
\DickPunkt(6,-10)
\DickPunkt(7,-10)
\DickPunkt(8,-10)
\DickPunkt(9,-10)
\DickPunkt(10,-10)
\DickPunkt(11,-10)
\DickPunkt(12,-10)
\DickPunkt(13,-10)
\DickPunkt(14,-10)
\Label\o{P}(8,-3)
\Label\o{\om}(0,0)
\Label\l{\bom}(0,-1)
\Label\ro{\kern7pt\om'}(13,-9)
\Label\r{\hskip2pt\hbox to 1cm{\leaders\hbox to
.3cm{\hss.\hss}\hfill}}(15,-10)
\Label\u{\hbox to 5.7cm{\leaders\hbox to
.3cm{\hss.\hss}\hfill}}(9,-10)
\hbox{\hskip5.1cm}
$$
\centerline{\eightpoint a. possible ``old" candidate cells, revisited\hskip1cm
b. possible ``new" candidate cells\hskip1cm}
\vskip10pt
\centerline{\eightpoint Figure 11}
}
\vskip10pt
Since $P$ is ``above" $\bP$, with no overlap of horizontal
pieces, we may weaken the above by saying that
``old" candidate cells are located in one of the following {\it two}
regions:
\roster
\item "(I)" In a row strictly below $\om'$; then the backward path of
the ``old" candidate cell must enter the row of $\om'$ weakly to the
left of $\om'$.
\item "(II)" In a row weakly above $\om'$,
and weakly to the left of $P$, but on $P$ only on its vertical
pieces, with their topmost cells excluded.
\endroster

In (HC2) the
transformations on the entries of $\bH$ are performed to obtain $H$. 
Because
of the definition of these transformations on $H$ in (HC2), the ``new"
candidate cells, i.e., those which arise if we jump into (HC*1) with
$(T,H)$ and distinguished cell $\om$, 
are related to the ``old" ones, i.e., those 
which arise if we jump into (HC*1) with
$(\bT,\bH)$ and distinguished cell $\bom$, 
as follows: Those ``old" candidate cells which
were located in a lower row than $\om'$ are 
``new" candidate cells as well. 
All the other ``old" candidate cells moved one unit to the left and
one unit up to form the corresponding ``new" candidate
cells. Therefore, ``new" candidate cells are located in one of the
following two regions (see Figure~11.b; as before, 
possible ``new" candidate cells are indicated by bold dots):
\roster
\item "(I)" In a row strictly below $\om'$; then the backward path of
the ``new" candidate cell must enter the row of $\om'$ weakly to the
left of $\om'$.
\item "(II)" In a row weakly above $\om'$,
and weakly to the left of $P$, but on $P$ only on its horizontal
pieces, with their rightmost cells excluded.
\endroster

Consequently, the backward path of a ``new" candidate cell is necessarily later
then $P$ in the total order of backward paths.

\smallskip
This completes the proof of the lemma.\quad \quad \qed
\enddemo

\subhead 4. Random generation of semistandard tableaux and of plane
partitions\endsubhead
The purpose of this section is to demonstrate that Algorithm~HC of
Section~2 can be used as a device for the random generation of
semistandard tableaux of a given shape and with bounded entries, and
for the random generation of plane partitions inside a given box.

Let $\la=(\la_1,\la_2,\dots,\la_r)$ be a
partition and $b$ an integer $\ge r$.
We claim that semistandard tableaux of shape $\la$ with entries
between $1$ and $b$ can be randomly generated as follows:

Choose a content tabloid $C$ of shape $\la$ at random. This is, of
course, easy, as we may choose the entry in each cell $\rh$ of $C$
between $1-c_\rh$ and $b$ independently and at random. Now apply
Algorithm~HC to $C$. This gives a pair $(T,H)$. Discard the hook
tabloid $H$, and define the output of this procedure to be the
semistandard tableau $T$. (Since we discard $H$ at the end, we 
actually do not have to worry about $H$ during the execution of
Algorithm~HC. I.e., the procedure may be streamlined
further by not assigning anything to $H$ in (HC0) and by disregarding
the operations on $H$ every time (HC2) is performed.)

$T$ is indeed a random semistandard tableaux of shape $\la$ with entries
between $1$ and $b$, because, as demonstrated before, Algorithm~HC is
a bijection between content tabloids $C$ and pairs $(T,H)$ where 
$T$ is a semistandard tableaux of shape $\la$ with entries
between $1$ and $b$, and where $H$ is a hook tabloid of shape $\la$. 
This implies in 
particular that any hook tabloid $H$ of shape $\la$ is equally likely as second
component of the pair $(T,H)$. Therefore, when discarded, any
semistandard tableau $T$ of shape $\la$ with entries
between $1$ and $b$ is equally likely as first component.

\medskip
Finally we turn to the random generation of plane partitions. Recall
that, given a partition $\la$, a {\it plane partition of shape
$\la$} is a filling
$P$ of the cells of $\la$
with integers such that the entries along rows and columns 
are weakly decreasing. We are interested here in plane partitions of
shape $(c^a)$ (this is an $a\times c$ rectangle, i.e., a Ferrers
diagram with $a$ rows of length $c$) with
entries between $0$ and $b$. Because plane partitions can also be
viewed geometrically as certain piles of unit cubes (cf\. 
\cite{\RobbAA}, for instance), these plane
partitions are also called {\it plane partitions inside an $a\times
b\times c$ box}.

Plane partitions inside an $a\times b\times c$ box are in bijection
with semistandard tableaux of shape $(c^a)$ with entries between $1$
and $a+b$. This is seen by replacing each entry, $e$ say, in the $i$-th
row of such a plane partition by $b-e+i$. Therefore Algorithm~HC, or
rather a proper deformation of it, can also be used for the random
generation of plane partitions inside a given box.
We give this deformation explicitly as Algorithm~PP below.

\bigskip
\noindent
{\bf Algorithm PP}. The input for the algorithm is a filling $F$
of the shape $(c^a)$ with the entry in cell $(i,j)$ between $i-a$ and
$b+j-1$, $i=1,2,\dots,a$, $j=1,2,\dots,c$.

\smallskip
(PP0) Set $P:=F$.
Call the last cell of $\la$ in the fixed total order of 
cells (i.e., the cell $(a,c)$) the {\it distinguished cell}, and
call the entry it contains {\it special}.
Continue with (PP1).

\smallskip
(PP1) Compare the special entry $s$ with its right neighbour,
$x$ say (if there is no right neighbour, then, by convention, 
we set $x=\infty$), and its bottom neighbour, $y$ say 
(if there is no bottom neighbour, then, by convention, we set $y=\infty$), 
see (4.1). If $s\ge x$ and $s\ge y$, then continue with (PP2).

If not, i.e\. if the special entry violates weak decrease along
rows or columns, then we have the following
situation,
$$\smatrix \format \sa\c\s\c\se\\
\hlinefor5\\
&s&&x&\\
\hlinefor5\\
&y&\\
\hlinefor3
\endsmatrix\ ,\tag4.1
$$
where at least one of $s<x$ and $s< y$ holds. If $x-1\ge y$ then do the move
$$\smatrix \format \sa\c\s\c\se\\
\hlinefor5\\
&x-1&&s&\\
\hlinefor5\\
&y&\\
\hlinefor3
\endsmatrix\ .\tag4.2
$$
If $x\le y$ then do the move
$$\smatrix \format \sa\c\s\c\se\\
\hlinefor5\\
&y&&x&\\
\hlinefor5\\
&s+1&\\
\hlinefor3
\endsmatrix\ .\tag4.3
$$
The ``new" special entry in (4.2) is $s$ again, the ``new" special entry 
in (4.3) is $s+1$.
Repeat (PP1). 

\smallskip
(PP2) Let $P$ be the array just
obtained. 
If $\om$ is the first cell of $(c^a)$ (i.e., the cell $(1,1)$), 
then stop. The output of the algorithm is $P$.

Otherwise, we choose as the
``new" distinguished cell the predecessor of the ``old" distinguished
cell in the fixed total
order of the cells, and as the ``new" special entry the entry it contains. 
Continue with (PP1).\quad \quad \qed

\medskip
There are already algorithms for the random generation of
semistandard tableaux of a given shape with bounded entries and of
plane partitions inside a given box in the literature. These are the
appropriate specializations of algorithms by Luby, Randall and
Sinclair \cite{\LuRSAA}, 
and by Propp and Wilson
\cite{\PropAD, \PrWiAA, \WilDAA}. The approaches of these
authors are however quite different. In \cite{\LuRSAA, \PropAD, \PrWiAA} 
they propose a method which is based on
a Markov chain approach called ``coupling from the past" (see also
\cite{\PrWiAB}). 
Wilson's algorithms \cite{\WilDAA} 
are based on a determinant algorithm originally
developped by Colbourn et al\. \cite{\CoMNAA}. (There are two
algorithms in \cite{\WilDAA}, one is based on Kasteleyn determinants,
the other on Gessel--Viennot determinants.)
It is then interesting to compare the performance of 
our algorithms with the algorithms
by Propp and Wilson. For example, the complexity of Algorithm~PP,
measured in terms of the number of modified jeu de taquin moves
(4.1)-(4.2) and (4.1)-(4.3) which has to be performed, is at worst
$O\big(\sum
_{i=a} ^{a+c-1}\binom i2\big)=O\big(c(3a^2+3ac+c^2)\big)$. (The average case
complexity is expected to be much better, probably quadratic in $a$
and $c$.) The complexity of
Wilson's algorithm based on Kasteleyn determinants \cite{\WilDAA,
Sec.~3} on the other hand is between $O\big((ab+bc+ca)^{3/2}\big)$
and $O\big((ab+bc+ca)^{5/2}\log^\al(ab+bc+ca)\big)$, depending on the
size of $a,b,c$, while the complexity of Wilson's algorithm based on
Gessel--Viennot determinants \cite{\WilDAA, Sec.~5} is
$O(b^{1.688}(ab+bc+ca))$. The complexity of the ``coupling from
the past" algorithms \cite{\LuRSAA, \PropAD, \PrWiAA} is
$O\big((a+b)^2(ab+bc+ca)\log(ab+bc+ca)\big)$ (this is partially
conjectural; see \cite{\WilDAC,
Sec.~5}).
So, the performance of our algorithm is (at least) in the range of the
others if $a,b,c$ are roughly of the same
magnitude, but should be better if one of the quantities $a,b,c$ is much
larger than the other two, because there is no dependence on $b$ in
the performance of our algorithm. 
It has to be emphasized, however, that 
our algorithm is rather specialized, and does not extend to cases
that are still covered by the algorithms by Luby, Randall, and Sinclair,
respectively by Propp and Wilson.

\subhead Appendix\endsubhead
The appendix contains a complete example for Algorithms~HC and HC* for
$\la=(4,3,3,2)$ and $b=7$, setting up a mapping between the two sides
of Figure~3. See the specific descriptions given in Examples~HC and HC*
of how to read the following
tables.

\newpage
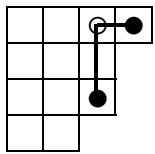
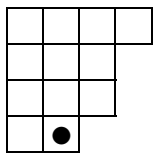
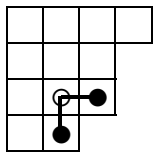
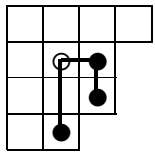
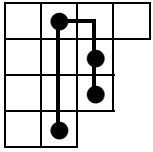
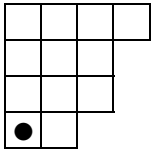
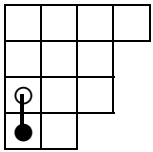
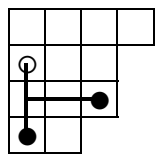
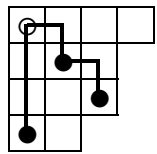
\eightpoint

\NoBlackBoxes
\Einheit.2cm
$$
\smatrix \format \l\s\l\\
\Ferrers
\Tableau
735{\GKreisum{--2}}/
732/
542/
46/
\
\Ferrers
\Tableau
0000/
000/
000/
00/
&&
\vtop{\hsize1.3cm
\vskip-1cm
$$\overset\hbox{\kern2.5pt(HC1)\kern2.5pt}\to\longrightarrow$$
$$\underset\hbox{\kern2.5pt(HC*2)\kern2.5pt}\to\longleftarrow$$
}
\Ferrers
\Tableau
735{\GBoxum{--2}}/
732/
542/
46/
\vtop{\hsize1.3cm
\vskip-1cm
$$\overset\hbox{\kern2.5pt(HC2)\kern2.5pt}\to\longrightarrow$$
$$\underset\hbox{\kern2.5pt(HC*1)\kern2.5pt}\to\longleftarrow$$
}
\\&&\\
\hskip46pt
&&
\\&&\\\hlinefor3\\&&\\
\Ferrers
\PfadDicke{1pt}
\Pfad(6,2),22\endPfad
\Tableau
735{\GDKreisum{--2}}/
732/
54{\Kreisum2}/
46/
\
\Ferrers
\Tableau
0000/
000/
000/
00/
&&
\vtop{\hsize1.3cm
\vskip-1cm
$$\overset\hbox{\kern2.5pt(HC1)\kern2.5pt}\to\longrightarrow$$
$$\underset\hbox{\kern2.5pt(HC*2)\kern2.5pt}\to\longleftarrow$$
}
\Ferrers
\Tableau
735{--2}/
732/
54{\Boxum2}/
46/
\vtop{\hsize1.3cm
\vskip-1cm
$$\overset\hbox{\kern2.5pt(HC2)\kern2.5pt}\to\longrightarrow$$
$$\underset\hbox{\kern2.5pt(HC*1)\kern2.5pt}\to\longleftarrow$$
}
\\&&\\
\Ferrers
\DickPunkt(7,3)
&&
\\&&\\\hlinefor3\\&&\\
\Ferrers
\PfadDicke{1pt}
\Pfad(4,-2),22112222\endPfad
\Tableau
735{--2}/
73{\Kreisum2}/
54{\DKreisum2}/
46/
\
\Ferrers
\Tableau
0000/
000/
000/
00/
&&
\vtop{\hsize1.3cm
\vskip-1cm
$$\overset\hbox{\kern2.5pt(HC1)\kern2.5pt}\to\longrightarrow$$
$$\underset\hbox{\kern2.5pt(HC*2)\kern2.5pt}\to\longleftarrow$$
}
\Ferrers
\Tableau
735{--2}/
73{\underbar2}/
542/
46/
\vtop{\hsize.9cm
\vskip-1cm
$$\overset\hbox{\kern2.5pt(2.4)\kern2.5pt}\to\longrightarrow$$
$$\underset\hbox{\kern2.5pt(2.8)\kern2.5pt}\to\longleftarrow$$
}
\Ferrers
\Tableau
735{--2}/
731/
54{\Boxum2}/
46/
\vtop{\hsize1.3cm
\vskip-1cm
$$\overset\hbox{\kern2.5pt(HC2)\kern2.5pt}\to\longrightarrow$$
$$\underset\hbox{\kern2.5pt(HC*1)\kern2.5pt}\to\longleftarrow$$
}
\\&&\\
\Ferrers
\DickPunkt(5,-1)
\hskip46pt
&&
\\&&\\\hlinefor3\\&&\\
\Ferrers
\PfadDicke{1pt}
\Pfad(4,-2),22221122\endPfad
\Tableau
73{\Kreisum5}{--2}/
731/
54{\DKreisum2}/
46/
\
\Ferrers
\Tableau
0000/
001/
000/
00/
&&
\vtop{\hsize1.3cm
\vskip-1cm
$$\overset\hbox{\kern2.5pt(HC1)\kern2.5pt}\to\longrightarrow$$
$$\underset\hbox{\kern2.5pt(HC*2)\kern2.5pt}\to\longleftarrow$$
}
\Ferrers
\Tableau
73{$\underline5$}{--2}/
731/
542/
46/
\vtop{\hsize.9cm
\vskip-1cm
$$\overset\hbox{\kern2.5pt(2.3)\kern2.5pt}\to\longrightarrow$$
$$\underset\hbox{\kern2.5pt(2.7)\kern2.5pt}\to\longleftarrow$$
}
\Ferrers
\Tableau
73{--1}{\Boxum5}/
731/
542/
46/
\vtop{\hsize1.3cm
\vskip-1cm
$$\overset\hbox{\kern2.5pt(HC2)\kern2.5pt}\to\longrightarrow$$
$$\underset\hbox{\kern2.5pt(HC*1)\kern2.5pt}\to\longleftarrow$$
}
\\&&\\
\Ferrers
\PfadDicke{1pt}
\DickPunkt(5,-1)
\Kreis(5,1)
\Pfad(5,-1),22\endPfad
\hskip46pt
&&
\\&&\\\hlinefor3\\&&\\
\endsmatrix
$$

\newpage
$$
\smatrix \format \l\s\l\\
\Ferrers
\PfadDicke{1pt}
\Pfad(4,-2),222222\endPfad
\Tableau
73{--1}{\DKreisum5}/
731/
542/
4{\Kreisum6}/
\
\Ferrers
\Tableau
00{--1}0/
001/
000/
00/
&&
\vtop{\hsize1.3cm
\vskip-1cm
$$\overset\hbox{\kern2.5pt(HC1)\kern2.5pt}\to\longrightarrow$$
$$\underset\hbox{\kern2.5pt(HC*2)\kern2.5pt}\to\longleftarrow$$
}
\Ferrers
\Tableau
73{--1}5/
731/
542/
4{\Boxum6}/
\vtop{\hsize1.3cm
\vskip-1cm
$$\overset\hbox{\kern2.5pt(HC2)\kern2.5pt}\to\longrightarrow$$
$$\underset\hbox{\kern2.5pt(HC*1)\kern2.5pt}\to\longleftarrow$$
}
\\&&\\
\Ferrers
\PfadDicke{1pt}
\Kreis(5,3)
\DickPunkt(5,-1)
\Pfad(5,-1),2222\endPfad
\DickPunkt(7,3)
\Pfad(5,3),11\endPfad
\hskip46pt
&&
\\&&\\\hlinefor3\\&&\\
\Ferrers
\PfadDicke{1pt}
\Pfad(2,-4),2211222222\endPfad
\Tableau
73{--1}5/
731/
5{\Kreisum{4}}2/
4{\DKreisum6}/
\
\Ferrers
\Tableau
00{--1}0/
001/
000/
00/
&&
\vtop{\hsize1.3cm
\vskip-1cm
$$\overset\hbox{\kern2.5pt(HC1)\kern2.5pt}\to\longrightarrow$$
$$\underset\hbox{\kern2.5pt(HC*2)\kern2.5pt}\to\longleftarrow$$
}
\Ferrers
\Tableau
73{--1}5/
731/
5{$\underline4$}2/
46/
\vtop{\hsize.9cm
\vskip-1cm
$$\overset\hbox{\kern2.5pt(2.3)\kern2.5pt}\to\longrightarrow$$
$$\underset\hbox{\kern2.5pt(2.7)\kern2.5pt}\to\longleftarrow$$
}
\Ferrers
\Tableau
73{--1}5/
731/
53{\Boxum4}/
46/
\vtop{\hsize1.3cm
\vskip-1cm
$$\overset\hbox{\kern2.5pt(HC2)\kern2.5pt}\to\longrightarrow$$
$$\underset\hbox{\kern2.5pt(HC*1)\kern2.5pt}\to\longleftarrow$$
}
\\&&\\
\Ferrers
\PfadDicke{1pt}
\DickPunkt(3,-3)
\hskip46pt
&&
\\&&\\\hlinefor3\\&&\\
\Ferrers
\PfadDicke{1pt}
\Pfad(2,-4),2222112222\endPfad
\Tableau
73{--1}5/
7{\Kreisum{3}}1/
53{\DKreisum4}/
46/
\
\Ferrers
\Tableau
00{--1}0/
001/
0{--1}0/
00/
&&
\vtop{\hsize1.3cm
\vskip-1cm
$$\overset\hbox{\kern2.5pt(HC1)\kern2.5pt}\to\longrightarrow$$
$$\underset\hbox{\kern2.5pt(HC*2)\kern2.5pt}\to\longleftarrow$$
}
\Ferrers
\Tableau
73{--1}5/
7{$\underline3$}1/
534/
46/
\vtop{\hsize.9cm
\vskip-1cm
$$\overset\hbox{\kern2.5pt(2.3)\kern2.5pt}\to\longrightarrow$$
$$\underset\hbox{\kern2.5pt(2.7)\kern2.5pt}\to\longleftarrow$$
}
\Ferrers
\Tableau
73{--1}5/
72{\Boxum3}/
534/
46/
\vtop{\hsize1.3cm
\vskip-1cm
$$\overset\hbox{\kern2.5pt(HC2)\kern2.5pt}\to\longrightarrow$$
$$\underset\hbox{\kern2.5pt(HC*1)\kern2.5pt}\to\longleftarrow$$
}
\\&&\\
\Ferrers
\PfadDicke{1pt}
\Kreis(3,-1)
\DickPunkt(3,-3)
\Pfad(3,-3),22\endPfad
\DickPunkt(5,-1)
\Pfad(3,-1),11\endPfad
\hskip46pt
&&
\\&&\\\hlinefor3\\&&\\
\Ferrers
\PfadDicke{1pt}
\Pfad(2,-4),2222221122\endPfad
\Tableau
7{\Kreisum{3}}{--1}5/
72{\DKreisum3}/
534/
46/
\
\Ferrers
\Tableau
00{--1}0/
0{--1}1/
0{--1}0/
00/
&&
\vtop{\hsize1.3cm
\vskip-1cm
$$\overset\hbox{\kern2.5pt(HC2)\kern2.5pt}\to\longrightarrow$$
$$\underset\hbox{\kern2.5pt(HC*1)\kern2.5pt}\to\longleftarrow$$
}
\Ferrers
\Tableau
7{$\underline3$}{--1}5/
723/
534/
46/
\vtop{\hsize.9cm
\vskip-1cm
$$\overset\hbox{\kern2.5pt(2.3)\kern2.5pt}\to\longrightarrow$$
$$\underset\hbox{\kern2.5pt(2.7)\kern2.5pt}\to\longleftarrow$$
}
\Ferrers
\Tableau
70{$\underline3$}5/
723/
534/
46/
\vtop{\hsize.9cm
\vskip-1cm
$$\overset\hbox{\kern2.5pt(2.4)\kern2.5pt}\to\longrightarrow$$
$$\underset\hbox{\kern2.5pt(2.8)\kern2.5pt}\to\longleftarrow$$
}
\Ferrers
\Tableau
7025/
72{\Boxum3}/
534/
46/
\vtop{\hsize1.3cm
\vskip-1cm
$$\overset\hbox{\kern2.5pt(HC2)\kern2.5pt}\to\longrightarrow$$
$$\underset\hbox{\kern2.5pt(HC*1)\kern2.5pt}\to\longleftarrow$$
}
\\&&\\
\Ferrers
\PfadDicke{1pt}
\Kreis(3,1)
\DickPunkt(5,1)
\Pfad(3,1),11\endPfad
\DickPunkt(5,-1)
\Pfad(5,-1),22\endPfad
\DickPunkt(3,-3)
\Pfad(3,-3),2222\endPfad
&&
\endsmatrix
$$

\newpage
$$
\smatrix \format \l\s\l\\
\Ferrers
\PfadDicke{1pt}
\Pfad(2,-4),22222222\endPfad
\Tableau
7025/
72{\DKreisum3}/
534/
{\Kreisum{4}}6/
\
\Ferrers
\Tableau
00{--1}0/
0{--1}1/
0{--1}0/
00/
&&
\vtop{\hsize1.3cm
\vskip-1cm
$$\overset\hbox{\kern2.5pt(HC1)\kern2.5pt}\to\longrightarrow$$
$$\underset\hbox{\kern2.5pt(HC*2)\kern2.5pt}\to\longleftarrow$$
}
\Ferrers
\Tableau
7025/
723/
534/
{\Boxum4}6/
\vtop{\hsize1.3cm
\vskip-1cm
$$\overset\hbox{\kern2.5pt(HC2)\kern2.5pt}\to\longrightarrow$$
$$\underset\hbox{\kern2.5pt(HC*1)\kern2.5pt}\to\longleftarrow$$
}
\\&&\\
\Ferrers
\PfadDicke{1pt}
\DickPunkt(3,3)
\DickPunkt(5,1)
\Pfad(3,3),11\endPfad
\Pfad(5,1),22\endPfad
\DickPunkt(5,-1)
\Pfad(5,-1),22\endPfad
\DickPunkt(3,-3)
\Pfad(3,-3),222222\endPfad
\hskip46pt
&&
\\&&\\\hlinefor3\\&&\\
\Ferrers
\PfadDicke{1pt}
\Pfad(0,-4),2211222222\endPfad
\Tableau
7025/
723/
{\Kreisum{5}}34/
{\DKreisum4}6/
\
\Ferrers
\Tableau
00{--1}0/
0{--1}1/
0{--1}0/
00/
&&
\vtop{\hsize1.3cm
\vskip-1cm
$$\overset\hbox{\kern2.5pt(HC1)\kern2.5pt}\to\longrightarrow$$
$$\underset\hbox{\kern2.5pt(HC*2)\kern2.5pt}\to\longleftarrow$$
}
\Ferrers
\Tableau
7025/
723/
{$\underline5$}34/
46/
\vtop{\hsize.9cm
\vskip-1cm
$$\overset\hbox{\kern2.5pt(2.4)\kern2.5pt}\to\longrightarrow$$
$$\underset\hbox{\kern2.5pt(2.8)\kern2.5pt}\to\longleftarrow$$
}
\Ferrers
\Tableau
7025/
723/
334/
{\Boxum5}6/
\vtop{\hsize1.3cm
\vskip-1cm
$$\overset\hbox{\kern2.5pt(HC2)\kern2.5pt}\to\longrightarrow$$
$$\underset\hbox{\kern2.5pt(HC*1)\kern2.5pt}\to\longleftarrow$$
}
\\&&\\
\Ferrers
\PfadDicke{1pt}
\DickPunkt(1,-3)
\hskip46pt
&&
\\&&\\\hlinefor3\\&&\\
\Ferrers
\PfadDicke{1pt}
\Pfad(0,-4),2222112222\endPfad
\Tableau
7025/
{\Kreisum{7}}23/
334/
{\DKreisum5}6/
\
\Ferrers
\Tableau
00{--1}0/
0{--1}1/
1{--1}0/
00/
&&
\vtop{\hsize1.3cm
\vskip-1cm
$$\overset\hbox{\kern2.5pt(HC1)\kern2.5pt}\to\longrightarrow$$
$$\underset\hbox{\kern2.5pt(HC*2)\kern2.5pt}\to\longleftarrow$$
}
\Ferrers
\Tableau
7025/
{$\underline7$}23/
334/
56/
\vtop{\hsize.9cm
\vskip-1cm
$$\overset\hbox{\kern2.5pt(2.4)\kern2.5pt}\to\longrightarrow$$
$$\underset\hbox{\kern2.5pt(2.8)\kern2.5pt}\to\longleftarrow$$
}
\Ferrers
\Tableau
7025/
223/
{$\underline7$}34/
56/
\vtop{\hsize.9cm
\vskip-1cm
$$\overset\hbox{\kern2.5pt(2.3)\kern2.5pt}\to\longrightarrow$$
$$\underset\hbox{\kern2.5pt(2.7)\kern2.5pt}\to\longleftarrow$$
}
\Ferrers
\Tableau
7025/
223/
4{$\underline7$}4/
56/
\vtop{\hsize.9cm
\vskip-1cm
$$\overset\hbox{\kern2.5pt(2.3)\kern2.5pt}\to\longrightarrow$$
$$\underset\hbox{\kern2.5pt(2.7)\kern2.5pt}\to\longleftarrow$$
}
\\&&\\
\Ferrers
\PfadDicke{1pt}
\Kreis(1,-1)
\DickPunkt(1,-3)
\Pfad(1,-3),22\endPfad
\hskip46pt
&&
\\&&\\\hlinefor3\\&&\\
&&
\vtop{\hsize.9cm
\vskip-1cm
$$\overset\hbox{\kern2.5pt(2.3)\kern2.5pt}\to\longrightarrow$$
$$\underset\hbox{\kern2.5pt(2.7)\kern2.5pt}\to\longleftarrow$$
}
\Ferrers
\Tableau
7025/
223/
45{\Boxum7}/
56/
\vtop{\hsize1.3cm
\vskip-1cm
$$\overset\hbox{\kern2.5pt(HC2)\kern2.5pt}\to\longrightarrow$$
$$\underset\hbox{\kern2.5pt(HC*1)\kern2.5pt}\to\longleftarrow$$
}
\\&&\\
\vbox to1cm{}
\endsmatrix
$$

\newpage
$$
\smatrix \format \l\s\l\\
\Ferrers
\PfadDicke{1pt}
\Pfad(0,-4),2222221122\endPfad
\Tableau
{\Kreisum{7}}025/
223/
45{\DKreisum7}/
56/
\
\Ferrers
\Tableau
00{--1}0/
2{--1}1/
{--2}{--1}0/
00/
&&
\vtop{\hsize1.3cm
\vskip-1cm
$$\overset\hbox{\kern2.5pt(HC1)\kern2.5pt}\to\longrightarrow$$
$$\underset\hbox{\kern2.5pt(HC*2)\kern2.5pt}\to\longleftarrow$$
}
\Ferrers
\Tableau
{$\underline7$}025/
223/
457/
56/
\vtop{\hsize.9cm
\vskip-1cm
$$\overset\hbox{\kern2.5pt(2.3)\kern2.5pt}\to\longrightarrow$$
$$\underset\hbox{\kern2.5pt(2.7)\kern2.5pt}\to\longleftarrow$$
}
\Ferrers
\Tableau
1{$\underline7$}25/
223/
457/
56/
\vtop{\hsize.9cm
\vskip-1cm
$$\overset\hbox{\kern2.5pt(2.4)\kern2.5pt}\to\longrightarrow$$
$$\underset\hbox{\kern2.5pt(2.8)\kern2.5pt}\to\longleftarrow$$
}
\Ferrers
\Tableau
1125/
2{$\underline7$}3/
457/
56/
\vtop{\hsize.9cm
\vskip-1cm
$$\overset\hbox{\kern2.5pt(2.3)\kern2.5pt}\to\longrightarrow$$
$$\underset\hbox{\kern2.5pt(2.7)\kern2.5pt}\to\longleftarrow$$
}
\\&&\\
\Ferrers
\PfadDicke{1pt}
\Kreis(1,1)
\DickPunkt(1,-3)
\Pfad(1,-3),2222\endPfad
\DickPunkt(5,-1)
\Pfad(1,-1),1111\endPfad
\hskip46pt
&&
\\&&\\\hlinefor3\\&&\\
&&
\vtop{\hsize.9cm
\vskip-1cm
$$\overset\hbox{\kern2.5pt(2.3)\kern2.5pt}\to\longrightarrow$$
$$\underset\hbox{\kern2.5pt(2.7)\kern2.5pt}\to\longleftarrow$$
}
\Ferrers
\Tableau
1125/
24{$\underline7$}/
457/
56/
\vtop{\hsize.9cm
\vskip-1cm
$$\overset\hbox{\kern2.5pt(2.4)\kern2.5pt}\to\longrightarrow$$
$$\underset\hbox{\kern2.5pt(2.8)\kern2.5pt}\to\longleftarrow$$
}
\Ferrers
\Tableau
1125/
246/
45{\Boxum7}/
56/
\vtop{\hsize1.3cm
\vskip-1cm
$$\overset\hbox{\kern2.5pt(HC2)\kern2.5pt}\to\longrightarrow$$
$$\underset\hbox{\kern2.5pt(HC*1)\kern2.5pt}\to\longleftarrow$$
}
\\&&\\
&&
\\&&\\\hlinefor3\\&&\\
\Ferrers
\Tableau
1125/
246/
45{\DKreisum7}/
56/
\
\Ferrers
\Tableau
30{--1}0/
{--1}{--1}1/
{--2}{--1}0/
00/
&&
\\&&\\
\Ferrers
\PfadDicke{1pt}
\Kreis(1,3)
\DickPunkt(1,-3)
\Pfad(1,-3),222222\endPfad
\DickPunkt(5,-1)
\DickPunkt(3,1)
\Pfad(1,3),11\endPfad
\Pfad(3,1),11\endPfad
\Pfad(3,1),22\endPfad
\Pfad(5,-1),22\endPfad
\hskip46pt
&&
\endsmatrix
$$


\newpage



\Refs

\ref\no \CoMNAA\by C. J. Colbourn, W. J. Myrvold and E. Neufeld \paper
Two algorithms for unranking arborescences\jour J. Algorithms
\vol 20\yr 1996\pages 268--281\endref

\ref\no \FrRTAA\by J. S. Frame, G. B. Robinson and R. M. Thrall \yr
1954 \paper The hook graphs of the symmetric group\jour
Canad\. J. Math\.\vol 6\pages 316--325\endref

\ref\no \GaMiAB\by A. M. Garsia and S. C. Milne \yr 1981 \paper Method
for constructing bijections for classical partition identities\jour
Proc\. Nat\. Acad\. Sci\. U.S.A.\vol 78\pages 2026--2028\endref

\ref\no \HiGrAA\by A. P. Hillman and R. M. Grassl \yr 1976 \paper
Reverse plane partitions and tableau hook numbers\jour J. Combin\.
Theory Ser.~A\vol 21\pages 216--221\endref

\ref\no \KratAY\by C.    Krattenthaler \yr 1998 \paper An involution
principle-free bijective proof of Stanley's hook-content formula\jour
Discrete Math\. Theoret\. Computer Science \toappear\vol \pages \endref

\ref\no \LuRSAA\by M. Luby, D. Randall and A. Sinclair\paper Markov
chain algorithms for planar lattice structures\paperinfo (extended
abstract)\inbook 36th Annual Symposium on Foundations of Computer
Science \yr 1995\pages 150--159\endref

\ref\no \MacMAA\by P. A. MacMahon \yr 1960 \book Combinatory Analysis
\publ Cambridge University Press, 1916; reprinted by Chelsea, New York
\endref

\ref\no \NoPSAA\by J. C. Novelli, I. M. Pak and A. V. Stoyanovsii 
\yr 1997 \paper A direct bijective proof of the hook-length formula\jour 
Discrete Math\. Theoret\. Computer Science\vol 1\pages 53--67\endref

\ref\no \PaStAA\by I. M. Pak and A. V. Stoyanovskii \yr 1992 \paper A
bijective proof of the hook-length formula and its analogues\jour
Funkt\. Anal\. Priloz\.\vol 26\rm , No.~3\pages 80--82\finalinfo
English translation in Funct\. Anal\. Appl\. {\bf 26} (1992),
216--218\endref

\ref\no \PropAD\by J.    Propp \yr 1997 \paper Generating random elements 
of  finite distributive lattices\jour Electron\. J. Combin\.\vol 
4 \rm(no.~2, ``The Wilf Festschrift")\pages \#R15, 12~pp\endref

\ref\no \PrWiAA\by J. Propp and D. B. Wilson\paper Exact sampling with
coupled Markov chains and applications to statistical mechanics\jour 
Random Structures Algorithms \vol 9\yr 1996\pages 223--252\endref

\ref\no \PrWiAB\by J.    Propp and D. B. Wilson \paper Coupling 
from the past: a user's guide\inbook Microsurveys in Discrete Probability
\eds D. Aldous, J. Propp\publ
DIMACS Ser\. Discrete Math\. Theoret\. Comput\. Sci\., vol.~41,
Amer\. Math\. Soc\., \publaddr Providence, R.I.\yr 1998
\pages 181--192\endref

\ref\no \ReWhAA\by J. B. Remmel and R. Whitney \yr 1983 \paper A
bijective proof of the hook formula for the number of column-strict
tableaux with bounded entries\jour Europ\. J. Combin\.\vol 4\pages
45--63\endref

\ref\no \RobbAA\by D. P. Robbins \yr 1991 \paper The story of 
$1,2,7,42,429,7436,\dots$\jour The Math\. Intelligencer\vol 13\pages 
12--19\endref

\ref\no \SagaAL\by B. E. Sagan \yr 1991 \book The symmetric group\publ
Wadsworth \& Brooks/Cole\publaddr Pacific Grove, California\endref

\ref\no \SchuAA\by M.-P. Sch\"utzenberger \yr 1977 \paper La
correspondance de Robinson\inbook Combinatoire et Repr\'esentation du
Groupe Sym\'etrique\publ Lecture Notes in Math\., vol.~579,
Springer--Verlag\publaddr Ber\-lin--Hei\-del\-berg--New York\pages
59--113\endref

\ref\no \StanAA\by R. P. Stanley \yr 1971 \paper Theory and
applications of plane partitions: Part~2\jour Stud\. Appl\. Math\vol
50\pages 259--279 \endref

\ref\no \WilDAA\by D. B. Wilson \yr 1997 \paper Determinant algorithms 
for random planar structures\inbook Proceedings of the Eighth Annual
ACM--SIAM Symposium on Discrete Algorithms\pages 258--267\endref

\ref\no \WilDAC\by D. B. Wilson \yr \paper Mixing times of lozenge
tiling and card shuffling Markov chains\paperinfo manuscript in
preparation\endref

\endRefs
\enddocument